\documentclass[a4paper,reqno, 11pt]{amsart}
\pdfoutput=1
\usepackage[utf8]{inputenc}
\usepackage[T1]{fontenc}
\usepackage{amssymb,mathtools}
\usepackage{mathrsfs}
\usepackage{dsfont}
\usepackage{bbm}

\usepackage[a4paper]{geometry}
\usepackage{lmodern}
\usepackage{color}
\usepackage{comment}
\usepackage[normalem]{ulem}
\usepackage{tikz}

\usepackage[all]{xy}
\usepackage{microtype}

\usepackage{enumitem}
\setlist[enumerate,1]{label=\textup{(\arabic*)}}

\usepackage[lite]{amsrefs}
\newcommand*{\MRref}[2]{ \href{http://www.ams.org/mathscinet-getitem?mr=#1}{MR \textbf{#1}}}

\renewcommand*{\PrintDOI}[1]{\href{http://dx.doi.org/\detokenize{#1}}{doi: \detokenize{#1}}}

\BibSpec{book}{%
  +{}  {\PrintPrimary}                {transition}
  +{,} { \textit}                     {title}
  +{.} { }                            {part}
  +{:} { \textit}                     {subtitle}
  +{,} { \PrintEdition}               {edition}
  +{}  { \PrintEditorsB}              {editor}
  +{,} { \PrintTranslatorsC}          {translator}
  +{,} { \PrintContributions}         {contribution}
  +{,} { }                            {series}
  +{,} { \voltext}                    {volume}
  +{,} { }                            {publisher}
  +{,} { }                            {organization}
  +{,} { }                            {address}
  +{,} { \PrintDateB}                 {date}
  +{,} { }                            {status}
  +{}  { \parenthesize}               {language}
  +{}  { \PrintTranslation}           {translation}
  +{;} { \PrintReprint}               {reprint}
  +{.} { }                            {note}
  +{.} {}                             {transition}
  +{} { \PrintDOI}                   {doi}
  +{} { available at \url}            {eprint}
  +{}  {\SentenceSpace \PrintReviews} {review}
}

\usepackage{hyperref}

\numberwithin{equation}{section}

\theoremstyle{plain}
\newtheorem{thm}[equation]{Theorem}
\newtheorem{cor}[equation]{Corollary}
\newtheorem{lem}[equation]{Lemma}
\newtheorem{prop}[equation]{Proposition}

\theoremstyle{definition}
\newtheorem{defn}[equation]{Definition}
\newtheorem{note}[equation]{Notation}

\theoremstyle{remark}
\newtheorem{rem}[equation]{Remark}

\newcommand{\NN}{\mathbb{N}}
\newcommand{\TT}{\mathbb{T}}
\newcommand{\CC}{\mathbb{C}}

\newcommand*{\nb}{\nobreakdash}
\newcommand*{\Star}{\(^*\)\nobreakdash-}

\newcommand{\Cst}{\mathrm{C}^*}

\newcommand{\idealin}{\mathrel{\triangleleft}} 
\newcommand*{\Bound}{\mathbb{B}}
\newcommand*{\Mat}{\mathbb{M}}
\newcommand*{\coloneqq}{\mathrel{\vcentcolon=}}
\newcommand{\Hilm}[1][E]{\mathcal{#1}}
\newcommand{\Toep}{\mathcal{T}}
\newcommand{\Toepr}{\mathcal{T}_\lambda}
\newcommand{\Toepu}{\mathcal{T}_u}
\newcommand{\CP}{\mathcal{O}}
\newcommand{\id}{\mathrm{id}}

\newcommand{\Fl}{\mathrm{C}^*_{\operatorname{rep}}(\mathcal{E})}

\newcommand*\xbar[1]{%
   \hbox{%
     \vbox{%
       \hrule height 0.5pt 
       \kern0.5ex
       \hbox{%
         \kern-0.1em
         \ensuremath{#1}%
         \kern-0.1em
       }%
     }%
   }%
} 

\DeclarePairedDelimiterX{\braket}[2]{\langle}{\rangle}{#1\,\delimsize\vert\,\mathopen{}#2}
\DeclarePairedDelimiterX{\BRAKET}[2]{\langle}{\rangle}{\!\delimsize\langle#1\,\delimsize\vert\,\mathopen{}#2\delimsize\rangle\!}

\DeclarePairedDelimiterX{\setgiven}[2]{\{}{\}}{#1\,{:}\,\mathopen{}#2}

\newcommand{\thmref}[1]{Theorem~\textup{\ref{#1}}}

\newcommand{\proref}[1]{Proposition~\textup{\ref{#1}}}
\newcommand{\lemref}[1]{Lemma~\textup{\ref{#1}}}
\newcommand{\corref}[1]{Corollary~\textup{\ref{#1}}}

\def\W{\mathcal W}

\def\1{\mathbbm 1}

\begin{document}
\title[$\Cst$-envelopes of tensor algebras]{$\Cst$-envelopes of tensor algebras of product systems}

\author{Camila F. Sehnem}

\address{School of Mathematics and Statistics, Victoria University of Wellington, P.O. Box 600, Wellington 6140, New Zealand.}

\thanks{Supported by the Marsden Fund of the Royal Society of New Zealand, grant No. \!\!18-VUW-056}

\thanks{\textit{Current address:} Department of Pure Mathematics, University of Waterloo, Waterloo, Ontario N2L 3G1, Canada.}

\thanks{\textit{Email address:} camila.sehnem@uwaterloo.ca}



\begin{abstract} Let $P$ be a submonoid of a group~$G$ and let $\Hilm=(\Hilm_p)_{p\in P}$ be a product system over~$P$ with coefficient $\Cst$\nb-algebra~$A$. We show that the following $\Cst$\nb-algebras are canonically isomorphic: the $\Cst$\nb-envelope of the tensor algebra~$\Toepr(\Hilm)^+$ of $\Hilm$; the reduced cross sectional $\Cst$\nb-algebra of the Fell bundle associated to the canonical coaction of~$G$ on the covariance algebra $A\times_{\Hilm}P$ of~$\Hilm$; and the $\Cst$\nb-envelope of the cosystem obtained by restricting the canonical gauge coaction on $\Toepr(\Hilm)$ to the tensor algebra. As a consequence, for every submonoid~$P$ of a group~$G$ and every product system $\Hilm=(\Hilm_p)_{p\in P}$ over~$P$, the $\Cst$\nb-envelope $\Cst_{\mathrm{env}}(\Toepr(\Hilm)^+)$ automatically carries a coaction of~$G$ that is compatible with the canonical gauge coaction on~$\Toepr(\Hilm)$. This answers a question left open by Dor-On, Kakariadis, Katsoulis, Laca and Li. We also analyse co-universal properties of $\Cst_{\mathrm{env}}(\Toepr(\Hilm)^+)$ with respect to injective gauge-compatible representations of~$\Hilm$. When $\Hilm=\CC^P$ is the canonical product  system over~$P$ with one-dimensional fibres, our main result implies that the boundary quotient $\partial\Toepr(P)$ is canonically isomorphic to the $\Cst$\nb-envelope of the closed non-selfadjoint subalgebra spanned by the canonical generating isometries of $\Toepr(P)$. Our results on co-universality imply that $\partial\Toepr(P)$ is a quotient of every nonzero $\Cst$\nb-algebra generated by a gauge-compatible isometric representation of $P$ that  in an appropriate sense respects the zero element of the semilattice of constructible right ideals of~$P$.

\end{abstract}
\maketitle

\section{Introduction}

Let $A$ be a $\Cst$\nb-algebra and let $\Hilm\colon A\leadsto A$ be a correspondence. The Toeplitz algebra~$\Toep_{\Hilm}$ of~$\Hilm$ is the $\Cst$\nb-algebra generated by the range of the canonical representation of~$\Hilm$ on $\Hilm^+=\bigoplus^\infty_{\substack{n\geq 0}}\Hilm^{\otimes n}$. Pimsner associated a $\Cst$\nb-algebra~$\CP_{\Hilm}$ to a faithful correspondence that unifies other constructions such as Cuntz algebras and crossed products by single automorphisms~\cite{Pimsner:Generalizing_Cuntz-Krieger}. Pimsner's $\Cst$\nb-algebra was generalised later by Katsura to not necessarily faithful correspondences \cite{Katsura:Cstar_correspondences}, and this construction of a $\Cst$\nb-algebra out of a correspondence is now known as a \emph{Cuntz--Pimsner algebra}. In a certain sense, the Cuntz--Pimsner algebra may be regarded as a crossed product of $A$ by the given correspondence.

The Cuntz--Pimsner algebra $\CP_{\Hilm}$ is a quotient of $\Toep_{\Hilm}$, and the corresponding quotient map is faithful on the copy of $A$. In fact, $\CP_{\Hilm}$ has the following special feature: the canonical gauge action of the unit circle on $\Toep_{\Hilm}$ induces a gauge action on~$\CP_{\Hilm}$ for which a representation of $\CP_{\Hilm}$ is faithful on the fixed-point algebra if and only if it is faithful on $A$ \cite[Proposition~6.3]{Katsura:Cstar_correspondences}. Although the definition of $\CP_{\Hilm}$ involves a notion of covariant representations,  $\CP_{\Hilm}$ can also be characterised as the quotient of $\Toep_{\Hilm}$ by its largest gauge-invariant ideal with trivial intersection with $A$ \cite[Proposition~7.14]{MR2413377}.

The copies of $\Hilm$ and $A$ in $\Toep_{\Hilm}$ generate a non-selfadjoint operator algebra $\Toep(\Hilm)^+$, called the \emph{tensor algebra} of $\Hilm$. The \emph{Shilov boundary} for $\Toep(\Hilm)^+$  is the largest ideal of $\Toep_{\Hilm}$ such that the corresponding quotient map is completely isometric on $\Toep(\Hilm)^+$. The quotient of $\Toep_{\Hilm}$ by this ideal is (isomorphic to) the \emph{$\Cst$\nb-envelope} of $\Toep(\Hilm)^+$, denoted by $\Cst_{\mathrm{env}}(\Toep(\Hilm)^+)$. The $\Cst$\nb-envelope of an operator system always exists by a result of Hamana \cite{Hamana}, following the seminal work of Arveson in \cite{Arveson}. The existence of the $\Cst$\nb-envelope for a non-unital operator algebra follows from the unital case and the work of Meyer on the unitization of an operator algebra~\cite{Meyer}.

While $\Cst_{\mathrm{env}}(\Toep(\Hilm)^+)$ is the smallest $\Cst$\nb-algebra generated by a completely isometric copy of $\Toep(\Hilm)^+$, the Cuntz--Pimsner algebra $\CP_{\Hilm}$ is the smallest $\Cst$\nb-algebra generated by an injective representation of $\Hilm$ admitting a gauge action of $\TT$ that is compatible with the canonical gauge action on $\Toep_{\Hilm}$. Muhly and Solel initiated in~\cite{Muhly-Solel:Tensor} the investigation on the relationship between these $\Cst$\nb-algebras. Under some assumptions on $\Hilm$, they proved that the canonical representation of $\Hilm$ in the $\Cst$\nb-envelope of $\Toep(\Hilm)^+$ induces an isomorphism  $\CP_{\Hilm}\cong\Cst_{\mathrm{env}}(\Toep(\Hilm)^+)$ \cite[Theorem~6.4]{Muhly-Solel:Tensor}. Katsoulis and Kribs established such an isomorphism  in the case of an arbitrary correspondence later in \cite{Katsoulis-Kribs}, following previous partial findings by Fowler, Muhly and Raeburn \cite[Theorem~5.3]{MR1986889}. Their result implies that  $\Cst_{\mathrm{env}}(\Toep(\Hilm)^+)$ automatically carries a gauge action of the unit circle~$\TT$ for which the quotient map from $\Toep_{\Hilm}$ is gauge-equivariant, and also connects notions of boundary quotients from the theories of selfadjoint and non-selfadjoint operator algebras.

Roughly speaking, a product system $\Hilm=(\Hilm_p)_{p\in P}$ over a monoid~$P$ as introduced by Fowler \cite{Fowler:Product_systems} is a semigroup of correspondences over the same $\Cst$\nb-algebra $A\coloneqq \Hilm_e$. A single correspondence naturally gives rise to a product system over $\NN$, and up to isomorphism every product system over $\NN$ is given by tensor powers of a single correspondence. If $P$ is left cancellative, there is a canonical representation of a product system $\Hilm=(\Hilm_p)_{p\in P}$ on the Hilbert $A$\nb-module $\Hilm^+=\bigoplus_{\substack{p\in P}}\Hilm_p$, called the Fock representation. This generates a $\Cst$\nb-algebra $\Toepr(\Hilm)$, called the Toeplitz algebra of~$\Hilm$. The non-selfadjoint closed subalgebra of $\Toepr(\Hilm)$ generated by the image of $\Hilm$ under the Fock representation is called the tensor algebra of $\Hilm$, and denoted by $\Toepr(\Hilm)^+$. As for a single correspondence, there are different notions of boundary quotients for a product system. In other words, there are different ways to associated a $\Cst$\nb-algebra to a product system that is the smallest $\Cst$\nb-algebra among those with a given property.

For $P$ a submonoid of a group~$G$, a notion of \emph{strongly covariant} representations of $\Hilm$ was introduced in \cite{SEHNEM2019558}, following previous work on compactly aligned product systems over positive cones of quasi-lattice orders by Sims and Yeend \cite{Sims-Yeend:Cstar_product_systems}. The \emph{covariance algebra} of $\Hilm$, denoted by $A\times_{\Hilm}P$, was then defined to be the universal $\Cst$\nb-algebra for strongly covariant representations of $\Hilm$. The $\Cst$\nb-algebra $A\times_{\Hilm}P$ does not depend on the group~$G$, and the canonical representation of $\Hilm$ in $A\times_{\Hilm}P$ is always injective. In addition, $A\times_{\Hilm}P$ satisfies a special property: it carries a canonical coaction of $G$ for which a representation of $A\times_{\Hilm}P$ is faithful on the fixed-point algebra if and only if it is faithful on the copy of $A$ \cite[Theorem~3.10]{SEHNEM2019558}. Thus the reduced analogue of  $A\times_{\Hilm}P$, namely the reduced cross sectional $\Cst$\nb-algebra of the Fell bundle associated to the canonical coaction on $A\times_{\Hilm}P$, is the smallest $\Cst$\nb-algebra generated by an injective strongly covariant representation of $\Hilm$ that admits a coaction of $G$ for which the quotient map from $\Toepr(\Hilm)$ is gauge-equivariant. Hence $A\times_{\Hilm}P$ may be viewed as a generalisation of a Cuntz--Pimsner algebra of a single correspondence.

 Dor-On, Kakariadis, Katsoulis, Laca and Li introduced in \cite{doron2020cenvelopes} the notion of a $\Cst$\nb-envelope for a cosystem. A cosystem consists of an operator algebra equipped with an appropriately defined coaction of a discrete group~$G$, and the $\Cst$\nb-envelope of a cosystem takes the coaction on the operator algebra into account. Hence the resulting $\Cst$\nb-algebra is the smallest $\Cst$\nb-algebra carrying a coaction of~$G$ that is generated by an equivariant completely isometric copy of the underlying operator algebra. For $P$ a submonoid of a group~$G$, the canonical gauge coaction of~$G$ on $\Toepr(\Hilm)$ restricts to a coaction on the tensor algebra $\Toepr(\Hilm)^+$ (see~\cite[Proposition~4.1]{doron2020cenvelopes}). The $\Cst$\nb-envelope of the corresponding cosystem has the $\Cst$\nb-envelope $\Cst_{\mathrm{env}}(\Toepr(\Hilm)^+)$ as a canonical quotient $\Cst$\nb-algebra.  However, the question of whether this quotient map is an isomorphism or, equivalently, the question of whether the $\Cst$\nb-envelope of~$\Toepr(\Hilm)^+$  automatically carries a coaction of $G$ for which the quotient map from~$\Toepr(\Hilm)$ is gauge-equivariant was left open in~\cite{doron2020cenvelopes}. By earlier results of Dor-On and Katsoulis, the quotient map was known to be an isomorphism in the case of compactly aligned product systems over positive cones of abelian lattice orders~\cite{MR4053621}. When $\Hilm$ is a compactly aligned product system over an arbitrary right LCM submonoid of a group, the $\Cst$\nb-envelope of the canonical cosystem associated to the tensor algebra~$\Toepr(\Hilm)^+$ was shown to be canonically isomorphic to the reduced analogue of the covariance algebra of~$\Hilm$ in \cite[Theorem~5.3]{doron2020cenvelopes}. But the relationship between these $\Cst$\nb-algebras remained unclear for more general product systems.

In this paper our main goal is to show that all of the above notions of boundary quotients for a product system $\Hilm=(\Hilm_p)_{p\in P}$ coincide. We were motivated by the recent analysis of the connection between these $\Cst$\nb-algebras in particular cases, see \cite{MR4053621, doron2020cenvelopes, kakariadis2021boundary, kakariadis2021couniversality}. Precisely, we show that the $\Cst$\nb-envelope of the tensor algebra $\Toepr(\Hilm)^+$ is canonically isomorphic to the reduced cross sectional $\Cst$\nb-algebra of the Fell bundle associated to the canonical coaction of~$G$ on the covariance algebra of~$\Hilm$, for $\Hilm=(\Hilm_p)_{p\in P}$ a product system over a submonoid~$P$ of~$G$. In particular, the $\Cst$\nb-envelope of $\Toepr(\Hilm)^+$ necessarily carries a coaction of~$G$, for which the quotient map from $\Toepr(\Hilm)$ is gauge-equivariant. Hence we also obtain that the $\Cst$\nb-envelopes of~$\Toepr(\Hilm)^+$ and of the cosystem arising from the canonical coaction on $\Toepr(\Hilm)^+$ are canonically isomorphic.

We point out that the original proof by Hamana of the existence of the $\Cst$\nb-envelope for an operator system is rather abstract \cite{Hamana}. In contrast, an important consequence of our main theorem is the identification of the $\Cst$\nb-envelope $\Cst_{\mathrm{env}}(\Toepr(\Hilm)^+)$ as the reduced analogue of the covariance algebra of $\Hilm$, which is concretely defined. This yields an explicit description of the Shilov boundary ideal for $\Toepr(\Hilm)^+$ in several situations, such as when the underlying group is known to be exact. See \corref{cor:Shilov}.

 We establish the existence of the isomorphisms mentioned above in \thmref{thm:main-result}. The proof builds on three main auxiliary results, which are proved in Sections \ref{sec:isometric} and \ref{sec:covariance}. First, we show in \corref{cor:gauge-iso} that the $\Cst$\nb-envelope of $\Toepr(\Hilm)^+$ is a quotient of every other quotient of $\Toepr(\Hilm)$ by a gauge-invariant ideal with trivial intersection with the coefficient algebra~$A$. In other words, any gauge-equivariant \Star homomorphism of $\Toepr(\Hilm)$ that is faithful on the coefficient algebra restricts to a complete isometry on $\Toepr(\Hilm)^+$. This implies that $\Cst_{\mathrm{env}}(\Toepr(\Hilm)^+)$ is a quotient of both the covariance algebra of $\Hilm$ and its reduced counterpart. The other two auxiliary results, namely \lemref{lem:iso-charac} and \proref{pro:char-expect}, allow us to conclude that the image under the quotient map of the spectral subspaces associated to the coaction on $\Toepr(\Hilm)$ yields a topological $G$\nb-grading for $\Cst_{\mathrm{env}}(\Toepr(\Hilm)^+)$. Thus the induced \Star homomorphism from the reduced analogue of the covariance algebra is an isomorphism.

After proving our main theorem in Section~\ref{sec:main-theorem}, we proceed to investigate co-universal properties for the $\Cst$\nb-envelope of $\Toepr(\Hilm)^+$ in Section~\ref{sec:co-universal}. We show that $\Cst_{\mathrm{env}}(\Toepr(\Hilm)^+)$ has the co-universal property for gauge-equivariant surjective \Star homomorphisms of the cross sectional $\Cst$\nb-algebra of the Fell bundle associated to the canonical coaction on $\Toepr(\Hilm)$ that are injective on $A$. See \corref{cor:gen-co-universal} for a precise statement. In particular this shows that the notion of $\Cst$\nb-envelopes of cosystems and the assumption of compact alignment in \cite[Theorem~4.9]{doron2020cenvelopes} are not needed to establish the existence of a $\Cst$\nb-algebra with this property. If~$\Hilm$ is faithful, meaning that the left action of~$A$ on each correspondence $\Hilm_p\colon A\leadsto A$ is injective, then we show that $\Cst_{\mathrm{env}}(\Toepr(\Hilm)^+)$ is a quotient of every $\Cst$\nb-algebra generated by an injective gauge-compatible representation of~$\Hilm$ that satisfies a certain condition, involving generating elements that correspond to the zero element of the semilattice of constructible right ideals of $P$ in a suitable sense (see \thmref{thm:zero-element}). This result may be particularly helpful when we wish to decide whether $\Cst_{\mathrm{env}}(\Toepr(\Hilm)^+)$ is a \Star homomorphic image of a given $\Cst$\nb-algebra. 

Although it is always helpful to have several descriptions for the same $\Cst$\nb-algebra at hand, we believe that this becomes very explicit in our analysis of co-universal properties in Section~\ref{sec:co-universal} in the case of the different descriptions for the $\Cst$\nb-envelope $\Cst_{\mathrm{env}}(\Toepr(\Hilm)^+)$ obtained from \thmref{thm:main-result}. The co-universal property of $\Cst_{\mathrm{env}}(\Toepr(\Hilm)^+)$ given in \corref{cor:gen-co-universal} follows as a consequence of its defining property as the smallest $\Cst$\nb-algebra generated by a completely isometric homomorphism of $\Toepr(\Hilm)^+$. On the other hand, the defining relations of the covariance algebra are our main tool in the proof of \thmref{thm:zero-element}, which gives a stronger version of co-universal property of $\Cst_{\mathrm{env}}(\Toepr(\Hilm)^+)$ when $\Hilm$ is faithful.

\section{Background}

In this section we review basic notions and results in the setting of non-selfadjoint operator algebras and of product systems that will be needed in the sequel, taking also the opportunity to establish our notation. We refer the reader to \cite{MR2111973, Vaulsen} for the theory of operator algebras and to \cite{Lance, MR1634408} for the theory of Hilbert $\Cst$\nb-modules. In this paper we will often use constructions of $\Cst$\nb-algebras associated to Fell bundles, and their relationship to topologically graded $\Cst$\nb-algebras in general. Our main reference for this part is \cite{Exel:Partial_dynamical}.

\subsection{Operator algebras and their \texorpdfstring{$\mathrm{C}^*$}{C*}-envelopes} An \emph{operator algebra} on a Hilbert space~$\Hilm[H]$ is a closed subalgebra of~$\Bound(\Hilm[H])$. If $A$ is an operator algebra on~$\Hilm[H]$, then for every $n\geq 1$ the algebra $\mathbb{M}_n(A)$ of $n\times n$ matrices with entries in $A$ is an operator algebra on the direct sum $\Hilm[H]^n\cong\CC^n\otimes\Hilm[H]$. We say that a homomorphism $\rho\colon A\to B$ between operator algebras is \emph{completely contractive} (resp. \emph{completely isometric}) if the induced homomorphism $\rho_n\colon \mathbb{M}_n(A)\to \mathbb{M}_n(B)$ is contractive (resp. isometric) for all~$n\geq 1$. 

A \emph{$\Cst$\nb-cover} of an operator algebra $A$ is a pair $(B,\rho)$, where $B$ is a $\Cst$\nb-algebra and $\rho\colon A\to B$ is a completely isometric homomorphism such that $\Cst(\rho(A))=B$. The \emph{$\Cst$\nb-envelope} of an operator algebra $A$ is a $\Cst$\nb-cover $(\Cst_{\mathrm{env}}(A),\iota)$ satisfying the following property: if $(B,\rho)$ is a $\Cst$\nb-cover of $A$, then there exists a (necessarily unique and surjective) \Star homomorphism $\pi\colon B\to \Cst_{\mathrm{env}}(A)$ such that $\pi\circ \rho=\iota$. The $\Cst$\nb-envelope of $A$ exists, and is unique up to an isomorphism that identifies $A$. Hamana established the existence of the $\Cst$\nb-envelope in the unital case in \cite[Theorem~4.1]{Hamana}, following the work of Arveson in \cite{Arveson}. The existence of the $\Cst$\nb-envelope in the non-unital case follows from the work of Meyer on the unitization of an operator algebra \cite{Meyer}.  See \cite[Proposition~4.3.5]{MR2111973} for further details.

Let $A$ be an operator algebra contained in a $\Cst$\nb-algebra $B$ and suppose that $\Cst(A)=B$. An ideal~$J$ in $B$ is called a \emph{boundary ideal} for $A$ if the quotient map $B\to B/J$ is completely isometric when restricted to $A$. We say that a boundary ideal $J$ of $B$ is the \emph{Shilov boundary} for $A$ if it contains every other boundary ideal for $A$. An ideal $J$ in $B$ is the Shilov boundary for $A$ if and only if the quotient $B/J$ is canonically isomorphic to the $\Cst$\nb-envelope of $A$.

\subsection{Product systems} Let $A$ be a $\Cst$\nb-algebra. A \emph{correspondence}~$\Hilm\colon A\leadsto A$ consists of a right Hilbert $A$\nb-module with a nondegenerate left action of~$A$ implemented by a \Star homomorphism $\varphi\colon A\rightarrow\Bound(\Hilm)$. We say that~$\Hilm$ is \emph{faithful} if the left action of~$A$ is injective. Notice that we have included nondegeneracy in our definition of a correspondence as we will be applying the main results of \cite{SEHNEM2019558}. By \cite[Remark 1.3]{Kwasniewski} this condition automatically holds for the underlying correspondences of a product system over a monoid with nontrivial group of units (see also \cite[Remark 2.2]{doron2020cenvelopes}).

Let~$P$ be a submonoid of a group $G$. We denote the unit element of $P$ by $e$. A \emph{product system} over $P$ of $A$\nb-correspondences consists of:
\begin{enumerate}
\item[(i)] a correspondence $\Hilm_p\colon A\leadsto A$ for each $p\in P\setminus\{e\}$;
\item[(ii)] correspondence isomorphisms $\mu_{p,q}\colon\Hilm_p\otimes_A\Hilm_q\overset{\cong}{\rightarrow}\Hilm_{pq}$, also called \emph{multiplication maps}, for all $p,q\in P\setminus\{e\}$;
\end{enumerate}

In addition, we let $\Hilm_e=A$ with the obvious structure of correspondence over~$A$. The multiplication maps $\mu_{e,p}$ and $\mu_{p,e}$ implement the left and right actions of~$A$ on~$\Hilm_p$, respectively. Thus~$\mu_{e,p}(a\otimes\xi_p)=\varphi_p(a)\xi_p$ and $\mu_{p,e}(\xi_p\otimes a)=\xi_pa$ for all~$a\in A$ and $\xi_p\in\Hilm_p$, where $\varphi_p\colon A\to\Bound(\Hilm_p)$ denotes the left action of~$A$ on $\Hilm_p$. 
 The multiplication maps must be associative, meaning that the following diagram commutes for all $p,q,r\in P$:
  \[
  \xymatrix{
    (\Hilm_p\otimes_A\Hilm_q)\otimes_A\Hilm_r  \ar@{->}[d]^{\mu_{p,q}\otimes1}
   \ar@{<->}[rr]& &     \Hilm_p\otimes_A(\Hilm_q\otimes_A\Hilm_r)
   \ar@{->}[rr]^{1\otimes\mu_{q,r}}&&
    \Hilm_p\otimes_A\Hilm_{qr} \ar@{->}[d]^{\mu_{p,qr}} \\
    \Hilm_{pq}\otimes_A\Hilm_r   \ar@{->}[rrrr]^{\mu_{pq,r}}&& &&
    \Hilm_{pqr}  
.  }
  \]
  
  A product system $\Hilm=(\Hilm_p)_{p\in P}$ will be called \emph{faithful} if $\varphi_p$ is injective for all $p\in P$. We refer to \cite{Fowler:Product_systems} for a more detailed account on product systems.

 \begin{defn}A \emph{representation} of a product system $\Hilm=(\Hilm_p)_{p\in P}$ in a $\Cst$\nb-algebra $B$ consists of linear maps~$\pi_p\colon\Hilm_p\rightarrow B$, for all $p\in P\setminus\{e\},$ and a \Star homomorphism $\pi_e\colon A\rightarrow B$, satisfying the following two axioms:
  \begin{enumerate}
  \item[(i)]$\pi_p(\xi)\pi_q(\eta)=\pi_{pq}(\xi\eta)$ for all $p,q\in P$, $\xi\in\Hilm_p$ and  $\eta\in\Hilm_q$;
  \item[(ii)] $\pi_p(\xi)^*\pi_p(\eta)=\pi_e(\braket{\xi}{\eta}),$ for all $p\in P$ and $\xi, \eta\in\Hilm_p$.  
  \end{enumerate}  
  \end{defn}

We will say that a representation $\pi=\{\pi_p\}_{p\in P}$ of a product system is \emph{injective} if the \Star homomorphism~$\pi_e$ is faithful. When $\pi=\{\pi_p\}_{p\in P}$ is injective, $\pi_p$ is automatically completely isometric for every $p\in P$.

\subsection{The Fock representation} There is a canonical injective representation associated to a product system $\Hilm=(\Hilm_p)_{p\in P}$, which we now describe. Let~$\Hilm^+$ be the right Hilbert $A$\nb-module given by the direct sum of all $\Hilm_p$'s. That is, $$\Hilm^+=\bigoplus_{\substack{p\in P}}\Hilm_p.$$ We call $\Hilm^+$ the \emph{Fock space} of $\Hilm$. For~$\xi\in \Hilm_p$, we define an operator $\psi^+_p(\xi)$ in $\Bound(\Hilm^+)$ by setting for each $\eta=\bigoplus_{\substack{s\in P}}\eta_s\in \Hilm^+$ $$\psi^+_p(\xi)(\eta)_s=\begin{cases}\mu_{p,p^{-1}s}(\xi\otimes\eta_{p^{-1}s})&\text{if } s\in pP,\\
0  &\text{otherwise}.
\end{cases}$$
We view~$\Hilm_{ps}$ as the correspondence $\Hilm_p\otimes_A\Hilm_s$ using the correspondence isomorphism~$\mu_{p,s}^{-1}$. In this way, $\psi_p^+(\xi)^*(\eta)_s$ is the image of~$\eta_{ps}$ in~$\Hilm_s$ under the operator defined on elements of the form $\mu_{p,s}(\zeta_p\otimes\zeta_s)$ by the formula $$\psi^+_p(\xi)^*(\mu_{p,s}(\zeta_p\otimes\zeta_s))=\varphi_s(\braket{\xi}{\zeta_p})\zeta_s.$$ So~$\psi_p^+(\xi)^*$ is the adjoint of~$\psi_p^+(\xi)$. This gives a representation $\psi^+=\{\psi^+_p\}_{p\in P}$ of~$\Hilm$ in~$\Bound(\Hilm^+)$ called the \emph{Fock representation} of~$\Hilm$. This representation is injective because the action of $A$ on the direct summand $A=\Hilm_e\subset\Hilm^+$ is simply left multiplication in $A$. See \cite{Fowler:Product_systems} for further details.

\begin{defn}\label{def:toeplitz-tensor-alg} We will refer to the $\Cst$\nb-algebra generated by the range of $\psi^+$ as the \emph{Toeplitz algebra} of $\Hilm$, and will denote it by $\Toepr(\Hilm)$. The \emph{tensor algebra} of $\Hilm$, denoted by~$\Toepr(\Hilm)^+,$ is the closed subalgebra of~$\Toepr(\Hilm)$ generated by the range of $\psi^+$. That is, $$\Toepr(\Hilm)^+=\overline{\mathrm{span}}\{\psi_p(\xi)\mid p\in P, \xi\in\Hilm_p\}.$$
\end{defn}

The universal $\Cst$\nb-algebra for representations of $\Hilm$ as introduced by Fowler in \cite[Proposition~2.8]{Fowler:Product_systems} will also be useful in this paper. We will denote this $\Cst$\nb-algebra by $\Fl$, and we let $\tilde{t}=\{\tilde{t}_p\}_{p\in P}$ denote the universal representation of $\Hilm$ in $\Fl$. Recall that the pair $(\Fl,\tilde{t})$ satisfies the following universal property: if $\pi=\{\pi_p\}_{p\in P}$ is a representation of $\Hilm$ in a $\Cst$\nb-algebra $B$, then there exists a unique \Star homomorphism $\tilde{\pi}\colon \Fl\to B$ such that $\tilde{\pi}\circ\tilde{t}_p=\pi_p$ for all $p\in P$.

\begin{rem} We have a few comments regarding the notation and terminology adopted in this paper:
\begin{enumerate}

\item We use the notation $\Fl$ for Fowler's universal $\Cst$\nb-algebra for representations of $\Hilm$, and we simply refer to it as Fowler's $\Cst$\nb-algebra. This
differs from the notation and terminology used in Fowler's original paper \cite{Fowler:Product_systems}, where this $\Cst$\nb-algebra was denoted by $\Toep_{\Hilm}$ and called the Toeplitz algebra of $\Hilm$; this notation and terminology were later adopted in other references, such as \cite{doron2020cenvelopes, SEHNEM2019558}. The reason for our choice of notation is that we believe that the symbols $\Toep_{\Hilm}$ and $\Toep(\Hilm)$ are reminiscent of a universal analogue of $\Toepr(\Hilm)$, in the sense that we should expect a canonical isomorphism between $\Toep(\Hilm)$ and $\Toepr(\Hilm)$ in nice cases such as, for example, when the enveloping group of $P$ is amenable. But such an isomorphism does not exist in general, even in simple cases such as when $P=\NN^2$ (see \cite{MR1386163}).

\item Here we reserve the name ``Toeplitz algebra'' for $\Toepr(\Hilm)$; this aligns with the terminology
arising from semigroup $\Cst$\nb-algebras and also with Pimsner's original definition of the
Toeplitz algebra of a single correspondence $\Hilm_1$ as the $\Cst$\nb-algebra generated by what is now known as the Fock representation of $\Hilm_1$ (see \cite{laca-sehnem, Pimsner:Generalizing_Cuntz-Krieger}). We point out that in \cite{doron2020cenvelopes} $\Toepr(\Hilm)$ is called the \emph{Fock algebra} of $\Hilm$ (see \cite[Definition 2.5]{doron2020cenvelopes}).

\end{enumerate}
\end{rem}

\subsection{The coaction on the Toeplitz algebra} Let~$G$ be a discrete group and let $\{u_g\mid g\in G\}$ be the canonical unitary generators of the full group $\Cst$\nb-algebra~$\Cst(G)$. Let $\delta_G$ be the \Star homomorphism $\Cst(G)\rightarrow \Cst(G)\otimes \Cst(G)$ given by~$\delta_G(u_g)=u_g\otimes u_g$. A \emph{(full) coaction} of $G$ on a $\Cst$\nb-algebra~$B$ is a nondegenerate and injective \Star homomorphism $\delta\colon B\rightarrow B\otimes \Cst(G)$ satisfying the coaction identity $$(\delta\otimes\id_{\Cst(G)})\delta=(\id_B\otimes\delta_G)\delta.$$ We refer to the triple $(B,G,\delta)$ as a coaction. The coaction $(B,G,\delta)$ is called \emph{nondegenerate} if $$\overline{\delta(B)(1\otimes\Cst(G))}=B\otimes\Cst(G).$$ It is said to be \emph{normal} if the \Star homomorphism $(\id_{B}\otimes \lambda)\circ\delta\colon B\to B\otimes \Cst_\lambda(G)$ is injective, where $\Cst_\lambda(G)$ is the reduced group $\Cst$\nb-algebra of~$G$ and $\lambda\colon \Cst(G)\to \Cst_\lambda(G)$, $u_g\mapsto \lambda_g$ is the left regular representation. See \cite{Quigg:Discrete_coactions_and_bundles} and also~\cite[Definition A.21]{Echterhoff-Kaliszewski-Quigg-Raeburn:Categorical}. Notice that it is not known whether all full coactions of discrete groups are automatically nondegenerate, see \cite{kaliszewski_quigg_2016}.

Let~$(B,G,\delta)$ be a nondegenerate coaction. The spectral subspace at $g\in G$ is $B_g\coloneqq\{b\in B\mid \delta(b)=b\otimes u_g\}$. We call~$B_e=B^\delta$ the \emph{fixed-point algebra} for the coaction~$\delta$ of~$G$ on~$B$. The collection of subspaces $\{B_g\}_{g\in G}$ forms a \emph{topological $G$\nb-grading} for $B$. That is, $\{B_g\}_{g\in G}$ is a collection of linearly independent closed subspaces of $B$ with $B_gB_h\subset B_{gh}$ and $B_g^*=B_{g^{-1}}$,  the sum $\sum B_g$ is dense in $B$, and there exists a conditional expectation $E^\delta\colon B\to B^\delta$ that vanishes on $B_g$ for $g\neq e$. The coaction $(B,G,\delta)$ is normal if and only if $E^\delta$ is faithful (see \cite[Lemma~1.4]{Quigg:Discrete_coactions_and_bundles}).  

As in \cite{doron2020cenvelopes}, we will need in the subsequent sections the normal coaction of~$G$ on the Toeplitz algebra.

\begin{prop}[\cite{doron2020cenvelopes}*{Proposition~4.1}] Let $P$ be a submonoid of a group $G$ and let $\Hilm=(\Hilm_p)_{p\in P}$ be a product system over~$P$ with coefficient $\Cst$\nb-algebra $A$. Let $\psi^+=\{\psi_p^+\}_{p\in P}$ be the Fock representation of~$\Hilm$. Then there is a nondegenerate normal coaction $$\overline{\delta}\colon \Toepr(\Hilm)\to\Toepr(\Hilm)\otimes\Cst(G)$$ that sends a generator $\psi^+_p(\xi)$ to $\psi^+_p(\xi)\otimes u_p$, for $p\in P$ and $\xi\in\Hilm_p$. Moreover, the spectral subspace~$\Toepr(\Hilm)_g$ at $g\in G$ is the closed linear span of elements of the form $$\psi_{p_1}^+(\xi_{p_1})\psi_{p_2}^+(\xi_{p_2})^*\ldots\psi_{p_{2k-1}}^+(\xi_{p_{2k-1}})\psi_{p_{2k}}^+(\xi_{p_{2k}})^*,$$ where $k\in \NN$, $p_1p_2^{-1}\ldots p_{2k-1}p_{2k}^{-1}=g$ and $\xi_{p_i}\in \Hilm_{p_i}$ for all $i\in\{1,2,\ldots,2k\}$.

\end{prop}

We will write $E_\lambda$ for the (faithful) conditional expectation $E^{\overline{\delta}}\colon \Toepr(\Hilm)\to \Toepr(\Hilm)_e$ associated to $\overline{\delta}$.
 
Fowler's $\Cst$\nb-algebra $\Fl$ also carries a canonical nondegenerate coaction of $G$. Indeed, there is a canonical representation of~$\Hilm$ in $\Fl\otimes \Cst(G)$ that sends~$\xi_p\in\Hilm_p$ to $\tilde{t}(\xi_p)\otimes u_p$. This induces a \Star homomorphism $\widetilde{\delta}\colon\Fl\rightarrow\Fl\otimes \Cst(G)$ by the universal property of~$\Fl$.  The triple $(\Fl, G, \widetilde{\delta})$ is a coaction, and as for $\Toepr(\Hilm)$ the spectral subspace $\Fl_g$ at $g\in G$ is the closed linear span of elements of the form $$\tilde{t}(\xi_{p_1})\tilde{t}(\xi_{p_2})^*\ldots\tilde{t}(\xi_{p_{2k-1}})\tilde{t}(\xi_{p_{2k}})^*,$$ where $k\in \NN$, $p_1p_2^{-1}\ldots p_{2k-1}p_{2k}^{-1}=g$ and $\xi_{p_i}\in \Hilm_{p_i}$ for all $i\in\{1,2,\ldots,2k\}$. See, for example, \cite[Lemma~2.2]{SEHNEM2019558}.

\section{Completely isometric homomorphisms of $\mathcal{T}_\lambda(\mathcal{E})^+$}\label{sec:isometric}

Our main goal in this section is to show that a gauge-equivariant quotient of $\Toepr(\Hilm)$ contains a completely isometric copy of $\Toepr(\Hilm)^+$ as long as the quotient map is faithful on~$A$.

Throughout this section we let $P$ be a submonoid of a group $G$. Let $\Hilm=(\Hilm_p)_{p\in P}$ be a product system over~$P$ and let $A\coloneqq \Hilm_e$ be its coefficient $\Cst$\nb-algebra.

\begin{prop}\label{prop:gauge-ccp} Let $\pi=\{\pi_p\}_{p\in P}$ be an injective representation of $\Hilm$ in a $\Cst$\nb-algebra~$B$ and suppose that the map that sends $\psi^+_p(\xi)$ to $\pi_p(\xi)$ for $p\in P$ and $\xi\in\Hilm_p$ induces a completely contractive homomorphism $\hat{\pi}\colon \Toepr(\Hilm)^+\to B$.   Suppose, in addition, that there exists a conditional expectation $E_\pi\colon \tilde{\pi}( \Fl)\to \tilde{\pi}(\Fl_e)$ such that the diagram
$$\xymatrix{
\Fl\ar@{->}[r]^{\tilde{\pi}} \ar@{->}[d]_{E^{\tilde{\delta}}}&
\tilde{\pi}(\Fl) \ar@{->}[d]_{E_\pi}\\
\Fl_e\ar@{->}[r]^{\tilde{\pi}_{}} &
\tilde{\pi}(\Fl_e)
    }$$ commutes, where $\tilde{\pi}\colon \Fl\to B$ is the \Star homomorphism obtained by universal property. Then $\hat{\pi}$ is completely isometric.

\begin{proof}  We begin by proving that $\hat{\pi}$ is isometric. To do so, it suffices to show that $$\|\sum_{\substack{s\in F_1}}\psi_s^+(\xi_s)\|=\|\sum_{\substack{s\in F_1}}\pi_s(\xi_s)\|,$$ where $F_1\subset P$ is a finite set and $\xi_s\in \Hilm_s$ for all $s\in F_1$. Let $F_2\subset P$ be another finite set and take an element $\eta=\sum_{\substack{r\in F_2}}\eta_r\in \Hilm^+$, where $\eta_r\in \Hilm_r$ for all $r\in F_2$. Consider the set $F\subset P$ given by $F\coloneqq\{sr\mid s\in F_1, r\in F_2\}$ and notice that $$\|\sum_{\substack{s\in F_1}}\psi_s^+(\xi_s)\eta\|^2=\|\braket{\sum_{\substack{s\in F_1}}\psi_s^+(\xi_s)\eta}{\sum_{\substack{s\in F_1}}\psi_s^+(\xi_s)\eta}\|=\|\sum_{\substack{p\in F}}\braket{\eta_p'}{\eta_p'}\|,$$ where $\eta_p'\coloneqq \sum_{\substack{sr=p}}\psi_s^+(\xi_s)\eta_r$ is simply the coordinate of  $\sum_{\substack{s\in F_1}}\psi_s^+(\xi_s)\eta$ at $p\in F$. Since $\pi$ is injective, we have that $\pi_e$ is isometric and so \begin{equation}\label{eq:vect-norm}
\|\sum_{\substack{p\in F}}\pi_e(\braket{\eta_p'}{\eta_p'})\|=\|\sum_{\substack{p\in F}}\braket{\eta_p'}{\eta_p'}\|=\|\sum_{\substack{s\in F_1}}\psi^+(\xi_s)\eta\|^2.\end{equation} Also, observe that 
\begin{equation*}\label{eq:mix-ineq}
\begin{aligned}
\sum_{\substack{p,q\in F}}\pi_p(\eta_p')^*\pi_q(\eta_q')&=\big(\sum_{\substack{r\in F_2}}\pi_r(\eta_r)\big)^*\big(\sum_{\substack{s,t\in F_1}}\pi_s(\xi_s)^*\pi_t(\xi_t)\big)\big(\sum_{r\in F_2}\pi_r(\eta_r)\big)\\&\leq\|\sum_{\substack{s\in F_1}}\pi_s(\xi_s)\|^2\sum_{\substack{r, r'\in F_2}}\pi_r(\eta_r)^*\pi_{r'}(\eta_{r'}).
\end{aligned}
\end{equation*}
 Applying the conditional expectation $E_\pi\colon \tilde{\pi}(\Fl)\to \tilde{\pi}(\Fl_e)$ to the resulting inequality and using that $E_\pi\circ \tilde{\pi}=\tilde{\pi}\circ E^{\tilde{\delta}}$, we deduce that \begin{equation}\label{eq:posi-ineq}\sum_{\substack{p\in F}}\pi_e(\braket{\eta_p'}{\eta_p'})\leq\|\sum_{\substack{s\in F_1}}\pi_s(\xi_s)\|^2\sum_{\substack{r\in F_2}}\pi_e(\braket{\eta_r}{\eta_r}).\end{equation} Notice that $$\|\sum_{\substack{r\in F_2}}\pi_e(\braket{\eta_r}{\eta_r})\|=\|\sum_{\substack{r\in F_2}}\braket{\eta_r}{\eta_r}\|=\|\eta\|^2,$$ and hence \eqref{eq:posi-ineq} together with \eqref{eq:vect-norm} yield the inequality \begin{equation*}\label{eq:ineq-norm}
\begin{aligned}
\|\sum_{\substack{s\in F_1}}\psi_s^+(\xi_s)\eta\|^2\leq \|\sum_{\substack{s\in F_1}}\pi_s(\xi_s)\|^2\|\eta\|^2.
\end{aligned}
\end{equation*} This implies that $\|\sum_{\substack{s\in F_1}}\psi_s^+(\xi_s)\|\leq \|\sum_{\substack{s\in F_1}}\pi_s(\xi_s)\|$ because the above inequality holds for every element $\eta\in\Hilm^+$ of finite support. Since $\hat{\pi}$ is contractive, we then conclude that $\hat{\pi}$ is isometric.

Next we show that $\hat{\pi}_n$ is isometric on $ \Mat_n(\Toepr(\Hilm)^+)$ also when $n\geq 2$. Let $b=(b_{i,j})\in \Mat_n(\Toepr(\Hilm)^+)$, where each entry $b_{i,j}\in\Toepr(\Hilm)^+$ is a finite sum $\sum_s\psi_s^+(\xi_s)$ with $\xi_s\in \Hilm_s$. Similarly, let $\eta=(\eta_1,\ldots,\eta_n)\in(\Hilm^+)^n$, where each $\eta_i$ is a finite sum of the form $\sum_{\substack{r}}\eta_{i,r}\in\Hilm^+$ with $\eta_{i,r}\in \Hilm_r$. We fix $i\in \{1,\ldots, n\}$ and for each $p\in P$ we let $\eta'_p$ be the coordinate of $\sum_{\substack{j=1}}^nb_{i,j}\eta_j\in\Hilm^+$ at $p\in P$. Then setting $F_i\coloneqq\{p\in P\mid \eta'_p\neq 0\}$, we see that $F_i$ is a finite set with $$\braket{\sum_{\substack{j=1}}^nb_{i,j}\eta_j}{\sum_{\substack{j=1}}^nb_{i,j}\eta_j}=\sum_{\substack{p\in F_i}}\braket{\eta_p'}{\eta_p'}.$$ Thus as in the case $n=1$ we have \begin{equation}\label{eq:n-vect-norm}\|b\eta\|^2=\|\sum_{\substack{i=1}}^n\braket{\sum_{\substack{j=1}}^nb_{i,j}\eta_j}{\sum_{\substack{j=1}}^nb_{i,j}\eta_j}\|=\|\sum_{\substack{i=1}}^n\sum_{\substack{p\in F_i}}\braket{\eta_p'}{\eta_p'}\|=\|\sum_{\substack{i=1}}^n\sum_{\substack{p\in F_i}}\pi_e(\braket{\eta_p'}{\eta_p'})\|.\end{equation} If we set $\pi(\eta_i)\coloneqq \sum_{r}\pi_r(\eta_{i,r})$ for $i=1,\ldots,n$ and let $\pi_{n,1}(\eta)\coloneqq (\pi(\eta_1),\ldots,\pi(\eta_n))$ regarded as an element of the direct sum $B^n$ with its $\ell_n^2$\nb-structure of Hilbert $B$\nb-module, we see that $$\sum_{\substack{p,q\in F_i}}\pi_{p}(\eta_p')^*\pi_q(\eta_q')=\big(\sum_{\substack{j=1}}^n\hat{\pi}(b_{i,j})\pi(\eta_j)\big)^*\big(\sum_{\substack{j=1}}^n\hat{\pi}(b_{i,j})\pi(\eta_j)\big),$$ and $$\braket{\pi_{n,1}(\eta)}{\pi_{n,1}(\eta)}=\sum_{\substack{i=1}}^n\pi(\eta_i)^*\pi(\eta_i)=\sum_{\substack{i=1}}^n\sum_{\substack{r,r'\in F_i}}\pi_r(\eta_{i,r})^*\pi_{r'}(\eta_{i,r'}).$$ Hence viewing $\hat{\pi}_n(b)$ as an adjointable operator on $B^n$ we conclude that \begin{equation*}
\begin{aligned}
\sum_{\substack{i=1}}^n\sum_{\substack{p,q\in F_i}}\pi_p(\eta_p')^*\pi_q(\eta_q')&=\braket{\hat{\pi}_n(b)\pi_{n,1}(\eta)}{\hat{\pi}_n(b)\pi_{n,1}(\eta)}\\&\leq \|\hat{\pi}_n(b)\|^2\braket{\pi_{n,1}(\eta)}{\pi_{n,1}(\eta)}\\&=\|\hat{\pi}_n(b)\|^2\sum_{\substack{i=1}}^n\sum_{\substack{r,r'\in F_i}}\pi_r(\eta_{i,r})^*\pi_{r'}(\eta_{i,r'}).
\end{aligned}
\end{equation*} Again applying $E_\pi$ to the resulting inequality and using that $E_\pi\circ \tilde{\pi}=\tilde{\pi}\circ E^{\tilde{\delta}}$, we obtain \begin{equation*}\label{eq:posi-ineqn}
\sum_{\substack{i=1}}^n\sum_{\substack{p\in F_i}}\pi_e(\braket{\eta_p'}{\eta_p'})\leq \|\hat{\pi}_n(b)\|^2\sum_{\substack{i=1}}^n\pi_e(\braket{\eta_i}{\eta_i}).
\end{equation*} Combining this with \eqref{eq:n-vect-norm} and observing that $\|\sum_{\substack{i=1}}^n\pi_e(\braket{\eta_i}{\eta_i})\|=\|\eta\|^2,$ we conclude that $$\|b\eta\|^2\leq \|\hat{\pi}_n(b)\|^2\|\eta\|^2.$$ This shows that $\|b\|\leq\|\hat{\pi}_n(b)\|$, and so the equality holds because $\hat{\pi}_n$ is contractive. Hence $\hat{\pi}_n$ is isometric for all $n\geq 1$, proving that $\hat{\pi}$ is completely isometric as asserted.
\end{proof}
\end{prop}

The following is a concrete situation in which \proref{prop:gauge-ccp} applies.

\begin{cor}\label{cor:gauge-iso} Let $\pi=\{\pi_p\}_{p\in P}$  be an injective representation of $\Hilm$ in  a $\Cst$\nb-algebra~$B$ and suppose that the map that sends $\psi^+_p(\xi)$ to $\pi_p(\xi)$ for $p\in P$ and $\xi\in\Hilm_p$ induces a \Star homomorphism $\hat{\pi}\colon\Toepr(\Hilm)\to B$.  Suppose, in addition, that there exists a conditional expectation $E_\pi\colon \hat{\pi}( \Toepr(\Hilm))\to \hat{\pi}(\Toepr(\Hilm)_e)$ such that the diagram
$$\xymatrix{
\Toepr(\Hilm)\ar@{->}[r]^{\hat{\pi}} \ar@{->}[d]_{E_\lambda}&
\hat{\pi}(\Toepr(\Hilm)) \ar@{->}[d]_{E_\pi}\\
\Toepr(\Hilm)_e\ar@{->}[r]^{\hat{\pi}_{}} &
\hat{\pi}(\Toepr(\Hilm)_e)
    }$$ commutes. Then the restriction of $\hat{\pi}$ to the tensor algebra $\Toepr(\Hilm)^+$ is completely isometric.
\end{cor}

\begin{rem} Notice that \corref{cor:gauge-iso} immediately implies that if $\hat{\pi}\colon \Toepr(\Hilm)\to B$ is a \Star homomorphism and $(B, G, \gamma)$ is a coaction for which $\hat{\pi}$ is $\overline{\delta}-\gamma$\nb-equivariant, in the sense that $\gamma\circ\hat{\pi}=(\hat{\pi}\otimes\id_{\Cst(G)})\circ\overline{\delta}$, then the restriction of $\hat{\pi}$ to the tensor algebra is completely isometric provided it is injective on~$A$.
\end{rem}

The next result will be an important tool in the proof that the $\Cst$\nb-envelope of $\Toepr(\Hilm)^+$ also carries a conditional  expectation as in Corollary~\ref{cor:gauge-iso}, where the representation of $\Hilm$ in $\Cst_{\mathrm{env}}(\Toepr(\Hilm)^+)$ in this case is simply the one induced by the inclusion $\Toepr(\Hilm)^+\hookrightarrow \Cst_{\mathrm{env}}(\Toepr(\Hilm)^+)$.
\begin{lem}\label{lem:iso-charac} Let $\pi=\{\pi_p\}_{p\in P}$ be a representation of $\Hilm$ in a $\Cst$\nb-algebra $B$ and suppose that the map that sends $\psi^+_p(\xi)$ to $\pi_p(\xi)$ for $p\in P$ and $\xi\in\Hilm_p$ induces a completely isometric homomorphism $\hat{\pi}\colon\Toepr(\Hilm)^+\to B$. Then for every $n\geq 1$, finite sets $F_1, F_2,\ldots, F_n\subset P$ and choice of elements $\xi_p\in\Hilm_p$ for $p\in F_i$, $i=1,\ldots,n$, we have that $$\|\sum_{\substack{i=1}}^n\sum_{\substack{p\in F_i}}\pi_e(\braket{\xi_p}{\xi_p})\|\leq \|\sum_{\substack{i=1}}^n\sum_{\substack{p,q\in F_i}}\pi_p(\xi_p)^*\pi_q(\xi_q)\|.$$
\begin{proof} Suppose that $\hat{\pi}$ is completely isometric. We consider first the case $n=1$. Let $F=F_1\subset P$ be a finite set and for each $p\in F$, take $\xi_p\in \Hilm_p$. Let $(u_\lambda)_{\lambda\in\Lambda}$ be an approximate identity for~$A$ and for each $\lambda \in\Lambda$, consider the element $\eta_\lambda\in\Hilm_e\subset\Hilm^+$ determined by $u_\lambda$. That is, the coordinate of $\eta_\lambda$ at $p\in P$ is $u_\lambda$ when $p=e$, and zero for $p\neq e$. Then $\|\eta_\lambda\|=\|u_\lambda\|=1$, and hence applying the operator $\sum_{\substack{p\in F}}\psi_p^+(\xi_p)\in\Bound(\Hilm^+)$ to $\eta_\lambda$ we deduce that $$\|\sum_{\substack{p\in F}}\braket{\xi_pu_\lambda}{\xi_pu_\lambda}\|=\|\big(\sum_{\substack{p\in F}}\psi_p^+(\xi_p)\big)\eta_\lambda\|^2\leq \|\sum_{\substack{p\in F}}\psi_p^+(\xi_p)\|^2.$$ Since $\hat{\pi}$ is isometric, the corresponding norm inequality holds in $B$ with $\pi$ in place of $\psi^+$, and so $$\|\sum_{\substack{p\in F}}\pi_e(\braket{\xi_pu_\lambda}{\xi_pu_\lambda})\|\leq \|\sum_{\substack{p\in F}}\pi_p(\xi_p)\|^2=\|\sum_{\substack{p,q\in F}}\pi_p(\xi_p)^*\pi_q(\xi_q)\|.$$ Taking the limit over $\lambda$ we obtain the desired inequality.

Now let $n\geq 2$, and for each $i=1,\ldots, n$, let $F_i\subset P$ be a finite set. For each $i\in\{1,\ldots,n\}$ and $p\in F_i$, choose $\xi_p\in \Hilm_p$. For each $i=1,\ldots,n$, set $b_i\coloneqq \sum_{\substack{p\in F_i}}\psi_p^+(\xi_p)$ and consider the matrix
 $$b\coloneqq  \left[\begin{array}{cccc}
     b_1& 0 &\cdots  & 0\\
          b_2 &0& \cdots &  0\\
        \vdots & \vdots & \vdots &  \vdots\\
    b_n & 0& \cdots &0\\
       \end{array}
\right]\in \Mat_n(\Toepr(\Hilm)^+).$$ That is, $b\in   \Mat_n(\Toepr(\Hilm)^+)$ is the matrix whose $i$th entry of the first column is $b_i$ and the remaining columns vanish. Again take an approximate identity $(u_\lambda)_{\lambda\in\Lambda}$ for $A$ and let $\eta_\lambda=u_\lambda\in\Hilm_e\subset\Hilm^+$ be as above. Let $\eta_{\lambda,n}\in(\Hilm^+)^n$ be given by $$\eta_{\lambda,n}\coloneqq (\eta_\lambda, \underbrace{0,\ldots, 0}_{n-1 \text{ times}}).$$  Then $\|\eta_{\lambda,n}\|=\|\eta_\lambda\|=1$, giving \begin{equation}\label{eq:mat-ineq}
\|b\,\eta_{\lambda,n}\|^2\leq \|b\|^2.
\end{equation}  Next we compute \begin{equation}\label{eq:vect-norm-2}\|b\,\eta_{\lambda,n}\|^2=\|\sum_{\substack{i=1}}^n\braket{b_i\eta_\lambda}{b_i\eta_\lambda}\|=\|\sum_{\substack{i=1}}^n\sum_{\substack{p\in F_i}}\braket{\xi_pu_\lambda}{\xi_pu_\lambda}\|,\end{equation} and observe that $\hat{\pi}_n(b)^*\hat{\pi}_n(b)$ is the matrix whose first entry of the first row is $$\sum_{\substack{i=1}}^n\hat{\pi}(b_i)^*\hat{\pi}(b_i)=\sum_{\substack{i=1}}^n\sum_{\substack{p,q\in F_i}}\pi_p(\xi_p)^*\pi_q(\xi_q)$$ and the remaining entries are all zero. Thus using that $\hat{\pi}$ is completely isometric, we can replace $\|b\|^2$ in \eqref{eq:mat-ineq} by $\|\hat{\pi}_n(b)\|^2$, and combining the resulting inequality with \eqref{eq:vect-norm-2} we obtain \begin{equation*}\|\sum_{\substack{i=1}}^n\sum_{\substack{p\in F_i}}\pi_e(\braket{\xi_pu_\lambda}{\xi_pu_\lambda})\|\leq\|\hat{\pi}_n(b)\|^2=\|\hat{\pi}_n(b)^*\hat{\pi}_n(b)\|=\|\sum_{\substack{i=1}}^n\sum_{\substack{p,q\in F_i}}\pi_p(\xi_p)^*\pi_q(\xi_q)\|.\end{equation*} Then the assertion follows by taking the limit over $\lambda$ in the inequality above.
\end{proof}
\end{lem}

\section{On covariance algebras and their gauge-equivariant quotients}\label{sec:covariance}

In this section we first present a few technical results concerning the norm of an element in the fixed-point algebra of the covariance algebra of a product system $\Hilm=(\Hilm_p)_{p\in P}$. Then we prove a result that will be used to show that the $\Cst$\nb-envelope of the tensor algebra of $\Hilm$ carries a coaction of $G$ for which the quotient map $\Toepr(\Hilm)\to \Cst_{\mathrm{env}}(\Toepr(\Hilm)^+)$ is gauge-equivariant. 

As usual we let $P$ be a submonoid of a group~$G$ and let $\Hilm=(\Hilm_p)_{p\in P}$ be a product system over~$P$ with coefficient $\Cst$\nb-algebra $A$.  We first recall the construction of the covariance algebra of $\Hilm=(\Hilm_p)_{p\in P}$ and its defining relations. We refer to \cite{SEHNEM2019558} for further details. Let $F\subset G$ be a finite set. We set $$K_F\coloneqq\underset{g\in F}{\bigcap}gP.$$ For each~$p\in P$ and each~$F\subset G$ finite, we define an ideal~$I_{p^{-1}(p\vee F)}\idealin A$ as follows. For 
$g\in F$, we let
$$I_{p^{-1}K_{\{p,g\}}}\coloneqq\begin{cases}\underset{r\in K_{\{p,g\}}}{\bigcap}\ker\varphi_{p^{-1}r}&\text{if } K_{\{p,g\}}\neq\emptyset\text{ and }p\not\in K_{\{p,g\}},\\
A& \text{otherwise. }
\end{cases}$$ 
We then set $$I_{p^{-1}(p\vee F)}\coloneqq\underset{g\in F}{\bigcap}I_{p^{-1}K_{\{p,g\}}}.$$  This yields a correspondence $\Hilm_{F}\colon A\leadsto A$ by setting \begin{equation*}\label{eq:finite_set_correspondence}\Hilm_{F}\coloneqq\bigoplus_{\substack{p\in P}}\Hilm_pI_{p^{-1}(p\vee F)}.
\end{equation*}

Let $\Hilm^+_F$ be the right Hilbert $A$\nb-module $\bigoplus_{\substack{g\in G}}\Hilm_{gF}$. For each $\xi\in \Hilm_p$, we define an operator $t_F^p(\xi)\in\Bound(\Hilm_F^+)$ so that it maps the direct summand~$\Hilm_{gF}$ into~$\Hilm_{pgF}$ for all~$g\in G$. Explicitly, $$t_F^p(\xi)(\eta_{r})\coloneqq\mu_{p,r}(\xi\otimes_A\eta_{r}),\quad\eta_{r}\in\Hilm_{r}I_{r^{-1}(r\vee gF)}.$$ This is well defined because $I_{r^{-1}(r\vee F)}=I_{(pr)^{-1}(pr\vee pF)}$ for each
$F\subset G$ finite and each~$p\in P$. Its adjoint $t_F^p(\xi)^*$ maps $\mu_{p,r}(\zeta_p\otimes \eta_r)$ to $\varphi_r(\braket{\xi}{\zeta_p})\eta_r$. Again this is well defined since~$I_{s^{-1}(s\vee F)}=I_{s^{-1}p(p^{-1}s\vee p^{-1}F)}$ for all~$s\in pP$.  Thus $t_F=\{t_F^p\}_{p\in P}$ is a representation of~$\Hilm$, and so it induces a \Star homomorphism $t_F\colon \Fl\rightarrow\Bound(\Hilm_F^+),$ still denoted by $t_F$ by abuse of language.

Let  $Q^F_g$ be the projection of $\Hilm_F^+$ onto the direct summand $\Hilm_{gF}$. Then $$t_F^p(\xi)Q^F_g=Q^F_{pg}t_F^p(\xi),\quad t_F^p(\xi)^*Q^F_g=Q^F_{p^{-1}g}t_F^p(\xi)^*$$ for all $p\in P$ and $\xi\in\Hilm_p$, which implies that $\Hilm_F$ is invariant under the image of $\Fl_e$ under $t_F$. If~$F_1\subset F_2$ are finite subsets of~$G$, then $$I_{p^{-1}(p\vee F_1)}\supset I_{p^{-1}(p\vee F_2)}$$ for all~$p\in P$ and hence~$\Hilm_{F_2}$ may be regarded as a closed submodule of~$\Hilm_{F_1}$. In particular, we have  $$\|Q^{F_2}_et_{F_2}(b)Q^{F_2}_e\|\leq \|Q^{F_1}_et_{F_1}(b)Q^{F_1}_e\|$$  for all $b\in \Fl_e$. So let~$F$ range in the directed set consisting of all finite subsets of~$G$ ordered by inclusion and define an ideal~$J_e\idealin\Fl_e$ by $$J_e\coloneqq\left\{b\in\Fl_e\middle|\;\lim_{\substack{F}}\|b\|_F=0\right\},$$ where $\|b\|_F\coloneqq \|Q_F^et_F(b)Q_F^e\|$. 

\begin{defn}[\cite{SEHNEM2019558}*{Definition~3.2}]\label{def:strong_covariance} A representation $\pi=\{\pi_p\}_{p\in P }$ of~$\Hilm$ in a $\Cst$\nb-algebra~$B$ is \emph{strongly covariant} if the induced \Star homomorphism $\tilde{\pi}\colon\Fl\to B$ vanishes on~$J_e$.
\end{defn}

Let $J_\infty\idealin\Fl$ be the ideal generated by~$J_e$. The \emph{covariance algebra} of $\Hilm$, denoted by $A\times_{\Hilm}P$, is the quotient $\Cst$\nb-algebra~$\Fl/J_\infty$.

We let $q\colon \Fl\to A\times_{\Hilm}P$ be the quotient map and let $j=\{j_p\}_{p\in P}$ be induced representation of $\Hilm$ in $A\times_{\Hilm}P$, that is, $j_p=q\circ\tilde{t}_p$ for all $p\in P$. Then the pair $(A\times_{\Hilm}P, j)$ has the following universal property: if $\pi=\{\pi_p\}_{p\in P}$ is a strongly covariant representation of $\Hilm$ in a $\Cst$\nb-algebra $B$, then there is a unique \Star homomorphism $\hat{\pi}\colon A\times_{\Hilm}P\to B$ such that $\hat{\pi}\circ j_p=\pi_p$ for all $p\in P$. It also follows that the group $G$ in question may be taken to be any group containing $P$ as a submonoid (see \cite[Theorem~3.10]{SEHNEM2019558}).

By \cite[Lemma~3.3]{SEHNEM2019558}, the ideal $J_\infty$ satisfies $$J_\infty=\overline{\bigoplus_{\substack{g\in G}}}( J_\infty\cap \Fl_g)=\overline{\bigoplus_{\substack{g\in G}}} \Fl_gJ_e,$$ and so the gauge coaction $(\Fl, G, \widetilde{\delta})$ gives rise to a canonical coaction $$\delta\colon A\times_{\Hilm}P\to (A\times_{\Hilm}P)\otimes\Cst(G)$$ with $\delta(j_p(\xi))=j_p(\xi)\otimes u_p$ for all $p\in P$ and $\xi\in \Hilm_p$. The spectral subspace $[A\times_{\Hilm}P]_g$ at $g\in G$ is canonically isomorphic to the quotient $\Fl_g/\Fl_gJ_e$. See  \cite[Lemma~3.4]{SEHNEM2019558} for further details. An important property of $A\times_{\Hilm}P$ is that a \Star homomorphism $\hat{\pi}\colon A\times_{\Hilm}P\to B$ is faithful on the fixed-point algebra $(A\times_{\Hilm}P)^\delta=[A\times_{\Hilm}P]_e$ for $\delta$ if and only if it is faithful on the coefficient algebra~$A$ (see \cite[Theorem~3.10]{SEHNEM2019558}(C3)).

In the next lemma we will give a more explicit description of the ideal $J_e$ and of the norm of an element in $(A\times_{\Hilm}P)^\delta$.

\begin{lem}\label{lem:limit-norms} Let $q\colon \Fl\to A\times_{\Hilm} P$ be the quotient map and for each finite set $F\subset G$, let $J_{e,F}$ be the kernel of the \Star homomorphism from $\Fl_e$ into $\Bound(\Hilm_F)$ given by $b\mapsto Q^F_et_F(b)Q^F_e$. Then
\begin{enumerate}
\item $J_e=\overline{\bigcup_{\substack{F}}J_{e,F}}$;

\item  $\|q(b)\|=\lim_{\substack{F\to\infty}}\|b\|_F$ for all $b\in\Fl_e$, where $F$ ranges over the finite subsets of~$G$ directed by inclusion.
\end{enumerate}
\begin{proof} We begin by proving (1). It follows from the definition of $J_e$ that $J_{e,F}\subset J_e$ for all $F\subset G$ finite. This gives the inclusion $\overline{\bigcup_{\substack{F}}J_{e,F}}\subset J_e$. In order to establish the reverse inclusion, let $b\in J_e$. Then $\lim_F\|b\|_F=0$, and thus given $\varepsilon>0$, there exists a finite set $F'\subset G$ such that $\|b\|_F< \varepsilon$ for all $F\subset G$ finite such that $F\supset F'$. Observe that $\|b\|_{F'}$ is precisely the norm of the image of $b$ in $\Fl_e/J_{e,F'}$ under the quotient map. Thus we can find $c\in J_{e,F'}$ such that $$\|b\|_{F'}=\inf\{\|b+d\|\mid d\in J_{e,F'}\}\leq \|b+c\|<\varepsilon.$$ This implies that $b$ lies in the closure of $\bigcup_{\substack{F}}J_{e,F}$, giving the inclusion $J_e\subset  \overline{\bigcup_{\substack{F}}J_{e,F}}$. 

Next we prove (2). Since $J_{e,F}\subset J_e$, the quotient map $\Fl_e\to \Fl/J_e=( A\times_{\Hilm} P)^\delta$ factors through $ \Fl_e/J_{e,F}$. Denoting by $q_F$ the quotient map $\Fl_e\to \Fl_e/J_{e,F}$, we deduce that for all $b\in \Fl_e$ $$\|q(b)\|\leq\|q_F(b)\|=\|b\|_F.$$ Because $F\subset G$ is an arbitrary finite set, it follows that $\|q(b)\|\leq \lim_{\substack{F\to \infty}}\|b\|_F$.

It remains to prove the inequality $\lim_{\substack{F\to \infty}}\|b\|_F\leq \|q(b)\|$. Given $\varepsilon>0$, let $c\in J_e$ be such that $\|q(b)\|\leq\|b+c\|<\|q(b)\|+\varepsilon$. Using (1), we may assume that $c\in J_{e,F'}$ for some finite set $F'\subset G$. It follows that $$\lim_{\substack{F\to \infty}}\|b\|_F\leq\|b\|_{F'}\leq \|b+c\|<\|q(b)\|+\varepsilon.$$ Since $\varepsilon>0$ is arbitrary, we conclude that $\lim_{\substack{F\to \infty}}\|b\|_F\leq \|q(b)\|$. This gives (2) and completes the proof of the lemma.
\end{proof}
\end{lem}

\begin{note}[see \cite{laca-sehnem}*{Section~2}] In what follows we will use the following notation:
\begin{itemize}
\item $\W(P)$ will denote the set of words in $P$ of even length.
\item For a word $\alpha=(p_1,p_2,\ldots, p_{2k-1},p_{2k})\in\W(P)$, we set $$\mathring{\alpha}\coloneqq p_1p_2^{-1}\ldots p_{2k-1}p_{2k}^{-1}\in G\quad\text{ and }\quad\tilde{\alpha}=(p_{2k},p_{2k-1},\ldots, p_{2},p_1)\in\W(P).$$

\item Given $\alpha\in \W(P)$, we set \begin{equation*}
\begin{aligned}
Q(\alpha)&\coloneqq\{ p_{2k}, p_{2k}p_{2k-1}^{-1}p_{2k-2},\ldots, \mathring{\tilde{\alpha}}p_1\},\quad \text{ and }\quad \\ K(\alpha)&\coloneqq K_{Q(\alpha)}=p_{2k}P\cap p_{2k}p_{2k-1}^{-1}p_{2k-2}P\cdots \mathring{\tilde{\alpha}}p_1P.
\end{aligned}
\end{equation*}
\end{itemize}

\end{note}

We refer to $Q(\alpha)$ as the \emph{iterated quotient set} of the word $\alpha$. We will say that the word $\alpha$ is \emph{neutral} if $\mathring{\alpha}=e$. The  set $K(\alpha)$ is a \emph{constructible right ideal} of $P$ in the sense of Li \cite{Li:Semigroup_amenability}.

\begin{rem} For the purpose of this paper it is more convenient to work with the generalised iterated \emph{right} quotient map $\W(P)\to G$ that associates the product of right quotients $\mathring{\alpha}=p_1p_2^{-1}\ldots p_{2k-1}p_{2k}^{-1}\in G$ to a word $\alpha=(p_1,\ldots, p_{2k})\in\W(P)$, rather than the iterated \emph{left} quotient map used in \cite[Section~2]{laca-sehnem}. This is because we consider elements of the form $$\psi_{p_1}^+(\xi_{p_1})\psi_{p_2}^+(\xi_{p_2})^*\ldots\psi_{p_{2k-1}}^+(\xi_{p_{2k-1}})\psi_{p_{2k}}^+(\xi_{p_{2k}})^*$$ as a canonical spanning set of $\Toepr(\Hilm)$. By doing so we are following~\cite{SEHNEM2019558} and~\cite{doron2020cenvelopes} and also avoid any further use of the nondegeneracy assumption on $\Hilm=(\Hilm_p)_{p\in P}$. As a result the meaning of $Q(\alpha)$ and $K(\alpha)$ here is slightly different from that in \cite[Section~2]{laca-sehnem}.
\end{rem}

\begin{lem}\label{lem:tech-lem} Let $j=\{j_p\}_{p\in P}$ be the canonical representation of $\Hilm$ in $A\times_{\Hilm}P$. Consider a word $\alpha=(p_1,p_2,\ldots, p_{2k-1},p_{2k})\in\W(P)$ and for each $i=1,\ldots, 2k$ let $\xi_{p_i}\in\Hilm_{p_i}$. Let $F\subset G$ be a finite set containing the iterated quotient set $Q(\alpha)$ of $\alpha$. 
\begin{enumerate}

\item If $r\in P$ and $r\not\in K(\alpha)$, then for all $\xi\in \Hilm_r I_{r^{-1}(r\vee F)}$ we have $$j_{p_1}(\xi_{p_1})j_{p_2}(\xi_{p_2})^*\ldots j_{p_{2k-1}}(\xi_{p_{2k-1}})j_{p_{2k}}(\xi_{p_{2k}})^*j_r(\xi)=0.$$

\item If  $r\in K(\alpha)$, then $\mathring{\alpha}r\in P$ and we have $$j_{p_1}(\xi_{p_1})j_{p_2}(\xi_{p_2})^*\ldots j_{p_{2k-1}}(\xi_{p_{2k-1}})j_{p_{2k}}(\xi_{p_{2k}})^*j_r(\Hilm_r I_{r^{-1}(r\vee F)})\subset j_{\mathring{\alpha}r}(\Hilm_{\mathring{\alpha}r}).$$ 
\end{enumerate} 
\begin{proof} Part (1) was shown in the proof of  \cite[Lemma~3.6]{SEHNEM2019558}, and so we briefly indicate how the proof goes. Suppose that $r\not\in K(\alpha)$ and let $\xi\in\Hilm_r I_{r^{-1}(r\vee F)}$. If $rP\cap K(\alpha)=\emptyset$,  then in $\Toepr(\Hilm)$ we have $$\psi^+_{p_1}(\xi_{p_1})\psi^+_{p_2}(\xi_{p_2})^*\ldots \psi^+_{p_{2k-1}}(\xi_{p_{2k-1}})\psi^+_{p_{2k}}(\xi_{p_{2k}})^*\psi^+_r(\xi)=0$$  because the restriction of $\psi^+_{p_1}(\xi_{p_1})\psi^+_{p_2}(\xi_{p_2})^*\ldots \psi^+_{p_{2k-1}}(\xi_{p_{2k-1}})\psi^+_{p_{2k}}(\xi_{p_{2k}})^*$ to a direct summand $\Hilm_s\subset \Hilm^+$ does not vanish only if $s\in K(\alpha)$. This implies that $$j_{p_1}(\xi_{p_1})j_{p_2}(\xi_{p_2})^*\ldots j_{p_{2k-1}}(\xi_{p_{2k-1}})j_{p_{2k}}(\xi_{p_{2k}})^*j_r(\xi)=0$$ in $A\times_{\Hilm}P$. Assume now that $rP\cap K(\alpha)\neq \emptyset$ but $r\not\in K(\alpha)$. Then there is $g\in Q(\alpha)$ such that $r\not\in gP$. Since $F\supset Q(\alpha)$, we also have that $g\in F$. It follows that $$\psi^+_{p_1}(\xi_{p_1})\psi^+_{p_2}(\xi_{p_2})^*\ldots \psi^+_{p_{2k-1}}(\xi_{p_{2k-1}})\psi^+_{p_{2k}}(\xi_{p_{2k}})^*\psi^+_r(\xi)=0$$ in $\Toepr(\Hilm)$ because if $s\in P$ is such that $rs\in K(\alpha)$, we have $rs\in rP\cap gP$ and so $I_{r^{-1}(r\vee F)}\subset \ker\varphi_s$. Thus the corresponding element in $A\times_{\Hilm}P$ also vanishes. This gives part (1). 

For part (2), notice that $K(\alpha)\subset\mathring{\tilde{\alpha}}p_1P\cap P$.  Thus if $r\in K(\alpha)$, we see that $\mathring{\alpha}r\in \mathring{\alpha}\mathring{\tilde{\alpha}}p_1P=p_1P\subset P$. So for $\xi\in \Hilm_r$, we deduce that $$j_{p_1}(\xi_{p_1})j_{p_2}(\xi_{p_2})^*\ldots j_{p_{2k-1}}(\xi_{p_{2k-1}})j_{p_{2k}}(\xi_{p_{2k}})^*j_r(\xi)\in j_{\mathring{\alpha}r}(\Hilm_{\mathring{\alpha}r})$$ because $j$ is a representation of $\Hilm$.
\end{proof}
\end{lem}

\begin{cor}\label{cor:normF} Let $q\colon \Fl\to A\times_{\Hilm}P$ be the quotient map and let $b\in\Fl_e$ be an element of the form $b=\sum_{\substack{i=1}}^{n}b_{i}$, where $n\geq 1$ and each $b_{i}$ in turn has the form $$b_{i}=\tilde{t}(\xi_{p_1})\tilde{t}(\xi_{p_2})^*\ldots \tilde{t}(\xi_{p_{2k-1}})\tilde{t}(\xi_{p_{2k}})^*$$ with $\alpha_i\coloneqq (p_1,\ldots, p_{2k})\in\W(P)$ a neutral word. Let $F\subset G$ be a finite set containing the iterated quotient set $Q(\alpha_i)$ for every $i=1,\ldots, n$. Then $$\|q(b)\|=\|b\|_F.$$
\begin{proof} It follows from part (2) of \lemref{lem:limit-norms} that $\|q(b)\|\leq\|b\|_F$. In order to show that $\|b\|_F\leq \|q(b)\|$, let $F'\subset P$ be a finite set and for each $r\in F'$ take $\xi_{r}\in  \Hilm_{r}I_{r^{-1}(r\vee F)}$.  Observe that because $F$ contains $Q(\alpha_i)$ for every $i=1,\ldots, n$, it follows from \lemref{lem:tech-lem} that for all $r\in F'$ and $i=1,\ldots, n$ we have that $q(b_i)j_{r}(\xi_{r})\in j_{r}(\Hilm_{r})$ in case $r\in K(\alpha_i)$, and  $q(b_i)j_{r}(\xi_{r})=0$ in case $r\not\in K(\alpha_i)$. Notice that if $r\in K(\alpha_i)$, then $q(b_i)j_{r}(\xi_{r})=j_r(\eta_r)$ with $\eta_r\coloneqq t_F(b_i)\xi_r\in\Hilm_rI_{r^{-1}(r\vee F)}$ since $j=\{j_p\}_{p\in P}$ is a representation of $\Hilm$. Now we have in $A\times_{\Hilm}P$ that \begin{equation*}
\begin{aligned}
\|\sum_{\substack{r\in F'}}j_{r}(\xi_{r})^*q(b)^*q(b)j_{r}(\xi_{r})\|\leq \|q(b)\|^2\|\sum_{\substack{r\in F'}}j_e(\braket{\xi_{r}}{\xi_{r}})\|.
\end{aligned}
\end{equation*} Hence if we let $\xi\coloneqq\sum_{\substack{r\in F'}}\xi_{r}\in\Hilm_{F}$, we see from the previous observation that the left-hand side above is precisely $\|t_{F}(b)\xi\|^2 $ (see also the proof of \cite[Lemma~3.6]{SEHNEM2019558}). Also, $\|\sum_{\substack{r\in F'}}j_e(\braket{\xi_{r}}{\xi_{r}})\|=\|\xi\|^2$. So we have $$\|t_{F}(b)\xi\|^2 \leq \|q(b)\|^2\|\xi\|^2,$$ which gives $\|b\|_{F}\leq \|q(b)\|$ because the set of elements of the form $\xi=\sum_{\substack{r\in F'}}\xi_{r}\in\Hilm_{F}$ with $F'\subset P$ finite is dense in $\Hilm_{F}$. This shows that $\|q(b)\|=\|b\|_F$ as asserted.
\end{proof}
\end{cor}

The next is the main result of this section.

\begin{prop}\label{pro:char-expect} Let $\pi=\{\pi_p\}_{p\in P}$ be an injective strongly covariant representation of~$\Hilm$ in a $\Cst$\nb-algebra~$B$ and let $\hat{\pi}\colon A\times_{\Hilm}P\to B$ be the induced \Star homomorphism. Then the following are equivalent:
\begin{enumerate}
\item There exists a conditional expectation $E_\pi\colon \hat{\pi}(A\times_{\Hilm} P)\to \hat{\pi}\big((A\times_{\Hilm} P)^\delta\big)$ such that the diagram $$\xymatrix{
A\times_{\Hilm}P\ar@{->}[r]^{\hat{\pi}} \ar@{->}[d]_{E^\delta}&
\hat{\pi}(A\times_{\Hilm} P) \ar@{->}[d]_{E_\pi}\\
(A\times_{\Hilm} P)^\delta\ar@{->}[r]^{\hat{\pi}} &
\hat{\pi}((A\times_{\Hilm} P)^\delta)
    }$$ commutes.

\item For every $n\geq 1$, finite sets $F_1, F_2,\ldots, F_n\subset P$ and elements $\xi_p\in\Hilm_p$ for $p\in F_i$, $i=1,\ldots,n$, we have that $$\|\sum_{\substack{i=1}}^n\sum_{\substack{p\in F_i}}\pi_e(\braket{\xi_p}{\xi_p})\|\leq \|\sum_{\substack{i=1}}^n\sum_{\substack{p,q\in F_i}}\pi_p(\xi_p)^*\pi_q(\xi_q)\|.$$ 

\end{enumerate}
\begin{proof} The implication (1)$\Rightarrow$(2) follows from evaluating $E_\pi$ at $ \sum_{\substack{i=1}}^n\sum_{\substack{p,q\in F_i}}\pi_p(\xi_p)^*\pi_q(\xi_q)$ and using that $E_\pi$ is contractive. Now assume that (2) holds. In order to establish the existence of $E_\pi$, it suffices to show that for every finite set $D\subset G$ and choice of elements $c_g\in [A\times_{\Hilm}P]_g$ for $g\in D$, we have $$\|\hat{\pi}(c_e)\|\leq \|\sum_{\substack{g\in D}}\hat{\pi}(c_g)\|.$$ This is so because the map $\sum\hat{\pi}(c_g)\mapsto \hat{\pi}(c_e)$ then extends by continuity to give a conditional expectation $E_\pi\colon \hat{\pi}(A\times_{\Hilm} P)\to \hat{\pi}\big((A\times_{\Hilm} P)^\delta\big)$ as in the statement. We may as well assume that each~$c_g$ has the form $c_g=q(b_g)$, where $b_g=\sum_{\substack{i=1}}^{n_g}b_{g,i}\in\Fl_g$, $n_g\geq 1$ and $b_{g,i}$ in turn has the form $$b_{g,i}=\tilde{t}(\xi_{p_1})\tilde{t}(\xi_{p_2})^*\ldots \tilde{t}(\xi_{p_{2k-1}})\tilde{t}(\xi_{p_{2k}})^*$$ with $\alpha_i\coloneqq (p_1,\ldots, p_{2k})\in\W(P)$ satisfying 
$\mathring{\alpha}_i=g$.

Let $c=\sum_{\substack{g\in D}}c_g\in A\times_{\Hilm}P$ and $b_g\in \Fl_g$ as above with $c_g=q(b_g)$. For each $g\in D$, consider the union of iterated quotient sets $F_g\coloneqq \cup_{\substack{i=1}}^{n_g}Q(\alpha_i)$. Put $F\coloneqq \bigcup_{\substack{g\in D}}F_g$. Fix $r\in P$ and let $\xi_r\in \Hilm_rI_{r^{-1}(r\vee F)}$. It follows from \lemref{lem:tech-lem} that $c_gj_r(\xi_r)=0$ if $gr\not\in P$ because then, in particular, $r\not\in K(\alpha_i)$ for all $i=1,\ldots, n_g$. In case $gr\in P$ we have $c_gj_r(\xi_r)\in j_{gr}(\Hilm_{gr})$ again by \lemref{lem:tech-lem}  because $F$ contains $Q(\alpha_i)$ for all $i=1,\ldots,n_g$. Consider the finite subset of $P$ given by $F_r\coloneqq \{gr\mid g\in D \text{ and }gr\in P\}$. For each $p\in F_r$ let $g\in D$ be such that $p=gr$ and let $\eta_p\in\Hilm_p$ be the vector satisfying $j_p(\eta_p)=c_gj_r(\xi_r).$ Then \begin{equation}\label{eq:1vector}
\begin{aligned}
\|\hat{\pi}(c)\|^2\pi_e(\braket{\xi_r}{\xi_r})&\geq \pi_r(\xi_r)^*\hat{\pi}(c)^*\hat{\pi}(c)\pi_r(\xi_r)\\&=\sum_{\substack{g,h\in D}}\pi_r(\xi_r)^*\hat{\pi}(c_g)^*\hat{\pi}(c_h)\pi_r(\xi_r)\\&=\sum_{\substack{p,q\in F_r}}\pi_p(\eta_p)^*\pi_q(\eta_q).
\end{aligned}
\end{equation} Thus applying (2) with $n=1$, $F_1=F_r$ and $\eta_p$ playing the role of $\xi_p$ for $p\in F_r$, we obtain the inequality $$\|\sum_{\substack{p\in F_r}}\pi_e(\braket{\eta_p}{\eta_p})\|\leq \|\hat{\pi}(c)\|^2\|\pi_e(\braket{\xi_r}{\xi_r})\|.$$ We compute \begin{equation}\label{eq:e-fibre}
\begin{aligned}
\sum_{\substack{p\in F_r}}\pi_e(\braket{\eta_p}{\eta_p})=\pi_r(\xi_r)^*\big(\sum_{\substack{g\in D}}\hat{\pi}(c_g)^*\hat{\pi}(c_g)\big)\pi_r(\xi_r)\geq \pi_r(\xi_r)^*\hat{\pi}(c_e)^*\hat{\pi}(c_e)\pi_r(\xi_r).
\end{aligned}
\end{equation}

Now let $F'\subset P$ be a finite set and for each $r\in F'$, take $\xi_{r}\in  \Hilm_{r}I_{r^{-1}(r\vee F)}$.  For each $r\in F'$, set $F_r=\{gr\mid g\in D\text{ and }gr\in P\}$. Then $F_r$ is a finite subset of $P$, and for each $p\in F_r$, take $g\in D$ such that $gr=p$. Let $\eta_p\in \Hilm_p$ be the vector such that $j(\eta_p)=c_gj(\xi_{r})$. Then as in \eqref{eq:1vector} we obtain 
 \begin{equation}\label{eq:nvector}
\begin{aligned}
\|\hat{\pi}(c)\|^2\sum_{\substack{r\in F'}}\pi_e(\braket{\xi_{r}}{\xi_{r}})&\geq\sum_{\substack{r\in F'}}\pi_r(\xi_{r})^*\hat{\pi}(c)^*\hat{\pi}(c)\pi_r(\xi_{r})\\&=\sum_{\substack{r\in F'}}\sum_{\substack{g,h\in D}}\pi_r(\xi_{r})^*\hat{\pi}(c_g)^*\hat{\pi}(c_h)\pi_r(\xi_{r})\\&=\sum_{\substack{r\in F'}}\sum_{\substack{p,q\in F_r}}\pi_p(\eta_p)^*\pi_q(\eta_q).
\end{aligned}
\end{equation} Thus applying (2) to the sum in the last line above with $n=|F'|$, the finite sets $F_r$ for $r\in F'$, and the elements $\eta_p\in \Hilm_p$ for $p\in F_r$ and using \eqref{eq:e-fibre}, we deduce that $$\|\sum_{\substack{r\in F'}}\pi_r(\xi_{r})^*\hat{\pi}(c_e)^*\hat{\pi}(c_e)\pi_r(\xi_{r})\|\leq\|\hat{\pi}(c)\|^2\|\sum_{\substack{r\in F'}}\pi_e(\braket{\xi_{r}}{\xi_{r}})\|.$$

As in the proof of Corollary~\ref{cor:normF}, let $\xi=\sum_{\substack{r\in F'}}\xi_{r}\in\Hilm_{F}$. Since $\pi=\{\pi_p\}_{p\in P}$ is injective, the left-hand side above is precisely $\|t_{F}(b_e)\xi\|^2 $, and also $\|\sum_{\substack{r\in F'}}\pi_e(\braket{\xi_{r}}{\xi_{r}})\|=\|\xi\|^2$. So we have $$\|t_{F}(b_e)\xi\|^2 \leq \|\hat{\pi}(c)\|^2\|\xi\|^2,$$ which implies that $\|b_e\|_{F}\leq \|\hat{\pi}(c)\|$. Combining this with Corollary~\ref{cor:normF} we obtain $$\|\hat{\pi}(c_e)\|=\|c_e\|=\|q(b_e)\|=\|b_e\|_F\leq \|\hat{\pi}(c)\|.$$ This implies the existence of $E_\pi$, completing the proof of the proposition.
\end{proof}
\end{prop}

\section{The main theorem}\label{sec:main-theorem}

We now combine the results from previous sections to show that the following three notions of boundary quotients associated to a product system $\Hilm=(\Hilm_p)_{p\in P}$ coincide: the $\Cst$\nb-envelope of the tensor algebra $\Toepr(\Hilm)^+$, in the sense of Arveson and Hamana; the $\Cst$\nb-envelope  of the canonical cosystem associated to $\Toepr(\Hilm)^+$ by Dor-On, Kakariadis, Katsoulis, Laca and Li \cite{doron2020cenvelopes}; and the reduced cross sectional $\Cst$\nb-algebra of the Fell bundle associated to the canonical coaction on the covariance algebra of $\Hilm$ as constructed in~\cite{SEHNEM2019558}.

 Notice that our approach only requires that $P$ be a submonoid of a group and our standard assumption that the left actions of~$A$ on the underlying correspondences are nondegenerate, since these are precisely the assumptions adopted in~\cite{SEHNEM2019558}. Thus as opposed to \cite{MR4053621, doron2020cenvelopes, kakariadis2021couniversality} we do not require that $P$ be a right LCM monoid, and even for $P$ a right LCM monoid, we do not assume $\Hilm$ to be compactly aligned. In case of compactly aligned product systems over right LCM monoids, we do not require the existence of a controlled map with abelian codomain as in \cite[Theorem~6.1]{kakariadis2021couniversality} (see \cite[Definition~5.1]{kakariadis2021couniversality}). Also, as a consequence of our main theorem, we deduce that the boundary quotient $\partial\Toepr(P)$ of the Toeplitz algebra $\Toepr(P)$ of a submonoid of a group is canonically isomorphic to the $\Cst$\nb-envelope of the closed non-selfadjoint subalgebra spanned by the canonical generating isometries of $\Toepr(P)$. This shows that the assumptions in Theorem~4.5 and Theorem~4.6 of~\cite{kakariadis2021boundary} are not necessary.

Before proving our main theorem, we briefly recall the definition of the $\Cst$\nb-envelope of a cosystem from \cite{doron2020cenvelopes}. Let $A$ be an operator algebra. A \emph{coaction} of a discrete group $G$ on $A$ is a completely isometric homomorphism $\gamma\colon A\to A\otimes \Cst(G)$ such that $\sum A_g$ is norm-dense in $A$, where $$A_g=\{a\in A\mid \gamma(a)= a\otimes u_g\}.$$ The triple $(A, G, \gamma)$ is called a \emph{cosystem}. See \cite[Definition~3.1]{doron2020cenvelopes} for further details. Notice that if $A$ is a $\Cst$\nb-algebra, $\gamma$ is automatically a nondegenerate coaction on $A$.

Let $(A, G, \gamma)$ and $(B, G, \gamma')$ be two cosystems. A completely contractive homomorphism $\phi\colon A\to B$ is said to be $\gamma-\gamma'$\nb-\emph{equivariant} or simply \emph{equivariant} if $\gamma'\circ\phi=(\phi\otimes\id_{\Cst(G)})\circ\gamma$. A triple $(B, \rho, \gamma')$ is a \emph{$\Cst$\nb-cover} for the cosystem $(A, G, \gamma)$ if $(B, G, \gamma')$ is a cosystem and the pair $(B,\rho)$ is a $\Cst$\nb-cover for $A$ with $\rho\colon A\to B$ being $\gamma$-$\gamma'$\nb-equivariant. See \cite[Definition~3.6]{doron2020cenvelopes}.

The \emph{$\Cst$\nb-envelope} of $(A, G, \gamma)$ is a $\Cst$\nb-cover for $(A, G, \gamma)$, denoted by $(\Cst_{\mathrm{env}}(A, G,\gamma), \iota_{\mathrm{env}},\gamma_{\mathrm{env}})$ or simply $\Cst_{\mathrm{env}}(A,G,\gamma)$, satisfying the following property: if $(B,\rho,\gamma')$ is a $\Cst$\nb-cover for $(A,G, \gamma)$, then there exists a $\gamma'-\gamma_{\mathrm{env}}$\nb-equivariant surjective \Star homomorphism $\phi\colon B\to \Cst_{\mathrm{env}}(A,G,\gamma)$ such that $\phi\circ\rho=\iota_{\mathrm{env}}$ \cite[Definition~3.7]{doron2020cenvelopes}. The $\Cst$\nb-envelope of a cosystem $(A,G,\gamma)$ always exists by \cite[Theorem~3.8]{doron2020cenvelopes}, and is unique up to a canonical isomorphism. 

Let $\Hilm=(\Hilm_p)_{p\in P}$ be a product system and let $\overline{\delta}\colon \Toepr(\Hilm)\to \Toepr(\Hilm)\otimes\Cst(G)$ be the canonical (normal) coaction of $G$ on $\Toepr(\Hilm)$. Let $\overline{\delta}^+$ denote the restriction of $\overline{\delta}$ to the tensor algebra $\Toepr(\Hilm)^+$. This gives rise to the cosystem $(\Toepr(\Hilm)^+, G, \overline{\delta}^+)$ as considered in \cite{doron2020cenvelopes, kakariadis2021couniversality, kakariadis2021boundary}.

\begin{thm}\label{thm:main-result} Let $P$ be a submonoid of a group $G$ and let $\Hilm=(\Hilm_p)_{p\in P}$ be a product system over $P$ with coefficient $\Cst$\nb-algebra~$A$. The following $\Cst$\nb-algebras associated to $\Hilm$ are canonically isomorphic:
\begin{enumerate}
\item the $\Cst$\nb-envelope $\Cst_{\mathrm{env}}(\Toepr(\Hilm)^+)$;

\item the reduced cross sectional $\Cst$\nb-algebra of the Fell bundle $([A\times_{\Hilm}P]_g)_{g\in G}$;

\item the $\Cst$\nb-envelope of the cosystem $(\Toepr(\Hilm)^+, G,\bar{\delta}^+)$.

\end{enumerate} 

\begin{proof} We begin by proving that the inclusion $\iota\colon \Toepr(\Hilm)^+\to \Cst_{\mathrm{env}}(\Toepr(\Hilm)^+)$ induces an isomorphism between the $\Cst$\nb-envelope  $\Cst_{\mathrm{env}}(\Toepr(\Hilm)^+)$ and the reduced cross sectional $\Cst$\nb-algebra $\Cst_r\big(([A\times_{\Hilm}P]_g)_{g\in G}\big)$. Consider the \Star homomorphism $\phi_\lambda\colon \Toepr(\Hilm)\to \Cst_r\big(([A\times_{\Hilm}P]_g)_{g\in G}\big)$ induced by the canonical representation of $\Hilm$ in $\Cst_r\big(([A\times_{\Hilm}P]_g)_{g\in G}\big)$ (for the existence of~$\phi_\lambda$ see, for example, the discussion after \cite[Proposition~5.4]{doron2020cenvelopes}). Then~$\phi_\lambda$ is injective on $A$, and the canonical conditional expectation $E_{\Lambda}$ of $\Cst_r\big(([A\times_{\Hilm}P]_g)_{g\in G}\big)$ onto the copy of $[A\times_{\Hilm}P]_e=(A\times_{\Hilm}P)^\delta$ satisfies $E_{\Lambda}\circ \phi_\lambda =\phi_{\lambda}\circ E_\lambda$. Thus it follows from \corref{cor:gauge-iso} that the restriction of $\phi_\lambda$ to the tensor algebra $\Toepr(\Hilm)^+$ is completely isometric. We deduce that there exists a \Star homomorphism $$\hat{\pi}\colon \Cst_r\big(([A\times_{\Hilm}P]_g)_{g\in G}\big)\to \Cst_{\mathrm{env}}(\Toepr(\Hilm)^+)$$ such that $\hat{\pi}\circ\phi_\lambda=\iota$ on $\Toepr(\Hilm)^+$. In particular, $\hat{\pi}$ is surjective. Since $\hat{\pi}$ is injective on $A$, it is also injective on the fixed-point algebra $[A\times_{\Hilm}P]_e= (A\times_{\Hilm}P)^\delta$ by \cite[Theorem~3.10]{SEHNEM2019558}. 

 The inclusion $\iota\colon \Toepr(\Hilm)^+\hookrightarrow \Cst_{\mathrm{env}}(\Toepr(\Hilm)^+)$ yields a canonical representation $\iota=\{\iota_p\}_{p\in P}$ of $\Hilm$ in $\Toepr(\Hilm)^+$, still denoted by $\iota$ by abuse of language. Since $\hat{\pi}\circ\phi_\lambda=\iota$ on $\Toepr(\Hilm)^+$, it follows that $\iota=\{\iota_p\}_{p\in P}$ is an injective strongly covariant representation of $\Hilm$. In addition, \lemref{lem:iso-charac} implies that $\iota=\{\iota_p\}_{p\in P}$ satisfies condition (2) of \proref{pro:char-expect} since $\iota\colon \Toepr(\Hilm)^+\hookrightarrow \Cst_{\mathrm{env}}(\Toepr(\Hilm)^+)$ is completely isometric. Hence there exists a conditional expectation $E_\iota\colon\Cst_{\mathrm{env}}(\Toepr(\Hilm)^+) \to\hat{\pi}((A\times_{\Hilm}P)^\delta)$ satisfying $E_\iota\circ\hat{\pi}=\hat{\pi}\circ E_\Lambda$. Because $E_{\Lambda}$ is faithful and $\hat{\pi}$ is faithful on $(A\times_{\Hilm}P)^\delta$, we deduce that $\hat{\pi}$ is an isomorphism by \cite[Theorem~19.5]{Exel:Partial_dynamical}.

Finally, let us establish the isomorphism between $\Cst_{\mathrm{env}}(\Toepr(\Hilm)^+)$ and the $\Cst$\nb-envelope of the cosystem $(\Toepr(\Hilm)^+, G,\bar{\delta}^+)$. It follows from the definition of $\Cst_{\mathrm{env}}(\Toepr(\Hilm)^+, G,\overline{\delta}^+)$ that there exists a surjective \Star homomorphism $$\phi\colon \Cst_{\mathrm{env}}(\Toepr(\Hilm)^+, G,\overline{\delta}^+)\to \Cst_{\mathrm{env}}(\Toepr(\Hilm)^+)$$  satisfying $\phi\circ\iota_{\mathrm{env}}=\iota$. To see that $\phi$ is an isomorphism, notice that Fell's absorption principle for Fell bundles, when combined with the canonical isomorphism $\Cst_{\mathrm{env}}(\Toepr(\Hilm)^+)\cong\Cst_r\big(([A\times_{\Hilm}P]_g)_{g\in G}\big)$ from the first part, gives that $\Cst_{\mathrm{env}}(\Toepr(\Hilm)^+)$ carries a (normal) coaction $$\delta_{\Lambda}\colon  \Cst_{\mathrm{env}}(\Toepr(\Hilm)^+)\to  \Cst_{\mathrm{env}}(\Toepr(\Hilm)^+)\otimes \Cst(G)$$ satisfying $\delta_{\Lambda}(\iota_p(\xi))= \iota_p(\xi)\otimes u_p$ for all $p\in P$ and $\xi\in\Hilm_p$ (see \cite[Proposition~18.24]{Exel:Partial_dynamical} and  \cite[Proposition~3.4]{doron2020cenvelopes}). Thus the inclusion $\iota\colon \Toepr(\Hilm)^+\hookrightarrow \Cst_{\mathrm{env}}(\Toepr(\Hilm)^+)$ is $\overline{\delta}^+-\delta_\Lambda$\nb-equivariant, which implies that $(\Cst_{\mathrm{env}}(\Toepr(\Hilm)^+),\iota,\delta_\Lambda)$ is a $\Cst$\nb-cover for the cosystem  $(\Toepr(\Hilm)^+, G,\bar{\delta}^+)$. Hence the defining property of  $\Cst_{\mathrm{env}}(\Toepr(\Hilm)^+, G,\bar{\delta}^+)$ yields the desired inverse for $\phi$, proving that $\phi$ is an isomorphism. 
\end{proof}
\end{thm}

\begin{rem} We observe that the defining relations of the covariance algebra form an essential part of the proof of  \proref{pro:char-expect}, and hence also of our proof that the $\Cst$\nb-envelope $\Cst_{\mathrm{env}}(\Toepr(\Hilm)^+)$ necessarily carries a coaction for which the inclusion of $\Toepr(\Hilm)^+$ is equivariant. In case of a more general operator algebra with a coaction, the question of whether the $\Cst$\nb-envelope of the corresponding cosystem coincides with the $\Cst$\nb-envelope of the underlying operator algebra remains open. More precisely, given a cosystem $(A, G, \gamma)$ that does not arise as the cosystem associated to the tensor algebra of a product system as above, it would be interesting to know whether the $\Cst$\nb-envelope of $A$ automatically carries a coaction of $G$ for which the inclusion $A\hookrightarrow \Cst_{\mathrm{env}}(A)$ is equivariant.
\end{rem}

We highlight in the next corollary an important consequence of \thmref{thm:main-result}. See also the discussion after Proposition~5.4 in~\cite{doron2020cenvelopes}.

\begin{cor}\label{cor:Shilov} Let $P$ be a submonoid of a group $G$ and let $\Hilm=(\Hilm_p)_{p\in P}$ be a product system over~$P$. Let $q_{\lambda}\colon \Fl\to\Toepr(\Hilm)$ be the \Star homomorphism induced by the Fock representation of~$\Hilm$. If the sequence $$0\xrightarrow{} J_\infty\xrightarrow{q_\lambda}\Toepr(\Hilm)\xrightarrow{\phi_\lambda}\Cst_r\big(([A\times_{\Hilm}P]_g)_{g\in G}\big)\xrightarrow{} 0$$ is exact (e.g. if $G$ is exact), then the Shilov boundary ideal for $\Toepr(\Hilm)^+$ is the ideal of $\Toepr(\Hilm)$ generated by  $$q_\lambda(J_e)=\{c\in\Toepr(\Hilm)_e\mid \lim_F\|c\|_F=0\},$$ where $F$ ranges over the finite subsets of $G$ ordered by inclusion and $\|c\|_F$ denotes the norm of the restriction of~$c$ to~$\Hilm_F$. 
\end{cor}

We consider now the particular case in which $\Hilm=\CC^P$ is the canonical product system over $P$ with one-dimensional fibres. The Fock representation of $\CC^P$ is simply the left regular representation $L\colon P\to \Bound(\ell^2(P))$ on $\ell^2(P)$. Thus $\Toepr(\CC^P)=\Toepr(P)$ is the Toeplitz algebra of $P$, and the corresponding tensor algebra of $\CC^P$, which we denote by $\Toepr(P)^+$, is then the closed non-selfadjoint subalgebra of $\Toepr(P)$ generated by the isometries $\{L_p\mid p\in P\}$. The fixed-point algebra $\Toepr(P)_e$ associated to the canonical coaction $\overline{\delta}$ of $G$ on $\Toepr(P)$ coincides with the diagonal $\Cst$\nb-subalgebra $D_r\subset \Toepr(P)$, and so $\Toepr(P)_e$ is a commutative $\Cst$\nb-algebra.  See, for example, \cite{CELY, laca-sehnem, Li:Semigroup_amenability} for further details on Toeplitz algebras of semigroups and other associated $\Cst$\nb-algebras.

There is a partial action $\gamma=(\{A_g\}_{g\in G},\{\gamma_g\})_{g\in G}$ of $G$ on $D_r$ such that $\Toepr(P)\cong D_r\rtimes_{\gamma, r}G$ canonically \cite[Theorem~5.6.41]{CELY} (see also \cite[Section~4]{laca-sehnem}). The Gelfand transform induces a partial action on $\mathrm{C}(\Omega_P)$ and an isomorphism $\Toepr(P)\cong \mathrm{C}(\Omega_P)\rtimes_r G$ in the natural way, where $\Omega_P$ denotes the spectrum of $D_r$. The \emph{boundary quotient} of $\Toepr(P)$, denoted by $\partial\Toepr(P)$, is the reduced partial crossed product $\mathrm{C}(\partial\Omega_P)\rtimes_r G$, where  $\partial\Omega_P$ is the smallest closed nonempty $\gamma$\nb-invariant subset of $\Omega_P$ (see \cite[Definition~5.7.9]{CELY}).

 It follows from \cite[Theorem~6.13]{laca-sehnem} (see also \cite[Theorem~3.14]{kakariadis2021boundary}) that $\partial\Toepr(P)$ is isomorphic to the reduced crossed sectional $\Cst$\nb-algebra of the Fell bundle associated to the canonical coaction of~$G$ on the covariance algebra of $\CC^P$, via an isomorphism that identifies the canonical generating isometries. That $\partial\Toepr(P)$ is also canonically isomorphic to the $\Cst$\nb-envelope of the cosystem $(\Toepr(P)^+, G,\overline{\delta}^+)$ was established in \cite[Theorem~4.4]{kakariadis2021boundary}. Although the authors proved in \cite{kakariadis2021boundary} that $\partial\Toepr(P)$ is canonically isomorphic to the $\Cst$\nb-envelope of $\Toepr(P)^+$ when either $\partial\Toepr(P)$ is simple or $P$ is an Ore monoid,  the question of whether $\partial\Toepr(P)$ is always canonically isomorphic to $\Cst_{\mathrm{env}}(\Toepr(P)^+)$ remained open even for group-embeddable right LCM monoids. We can now answer this question in the affirmative for arbitrary submonoids of groups as an immediate application of~\thmref{thm:main-result}.

\begin{cor}\label{cor:boundary-quotient} Let $P$ be a submonoid of a group. Then the boundary quotient $\partial\Toepr(P)$ is isomorphic to the $\Cst$\nb-envelope of $\Toepr(P)^+=\overline{\mathrm{span}}\{L_p\mid p\in P\}$, with an isomorphism that identifies the canonical generating isometries.
\end{cor}

\section{On co-universal properties for the \texorpdfstring{$\mathrm{C}^*$}{C*}-envelope $\mathrm{C}^*_{\mathrm{env}}(\mathcal{T}_\lambda(\mathcal{E})^+)$}\label{sec:co-universal}

For a compactly aligned product system $\Hilm=(\Hilm_p)_{p\in P}$ over a positive cone of a quasi-lattice order $(G,P)$, Carlsen, Larsen, Sims and Vittadello asked in \cite{Carlsen-Larsen-Sims-Vittadello:Co-universal} for the existence of a $\Cst$\nb-algebra satisfying a certain co-universal property with respect to injective gauge-compatible Nica covariant representations of $\Hilm$. In \cite[Theorem~4.1]{Carlsen-Larsen-Sims-Vittadello:Co-universal} they proved that under extra assumptions on $\Hilm$, the reduced analogue of the Cuntz--Nica--Pimsner algebra as defined by Sims and Yeend in \cite{Sims-Yeend:Cstar_product_systems}, denoted by $\mathcal{N}\CP^r_{\Hilm}$, satisfies the desired properties. That is, $\mathcal{N}\CP^r_{\Hilm}$ carries a coaction $\nu$ of $G$ for which the quotient map from the Nica--Toeplitz algebra $\mathcal{N}\Toep(\Hilm)$ onto $\mathcal{N}\CP^r_{\Hilm}$ is gauge-equivariant, and if $(B,G, \gamma)$ is a coaction and $B$ is generated as a $\Cst$\nb-algebra by an injective Nica covariant representation of $\Hilm$ that is gauge-compatible with $\gamma$ (see Definition~\ref{def:gauge-comp} below), then there exists a surjective \Star homomorphism $\rho\colon B\to \mathcal{N}\CP^r_{\Hilm}$ that identifies the corresponding copies of $\Hilm$. More recently,  in \cite[Theorem~4.9]{doron2020cenvelopes}  (see also \cite{MR4053621}) the authors proved that for $\Hilm$ a compactly aligned product system over a right LCM submonoid of $G$, the $\Cst$\nb-envelope of the cosystem $(\Toepr(\Hilm)^+, G,\overline{\delta}^+)$ always satisfies the co-universal property for injective gauge-compatible Nica covariant representations.

Motivated by the results mentioned above, we examine next co-universal properties for the $\Cst$\nb-envelope $\Cst_{\mathrm{env}}(\Toepr(\Hilm)^+)$ for product systems over arbitrary submonoids of groups. Notice that if $(\Toepr(\Hilm)_g)_{g\in G}$ is the Fell bundle associated to the coaction~$\bar{\delta}$ of~$G$ on~$\Toepr(\Hilm)$, then the cross sectional $\Cst$\nb-algebra $\Cst((\Toepr(\Hilm)_g)_{g\in G})$ is canonically isomorphic to the Nica--Toeplitz algebra of $\Hilm$ when $\Hilm$ is a compactly aligned product system over a right LCM submonoid of a group (see \cite[Proposition~4.3]{doron2020cenvelopes}).

\begin{defn}[see \cite{Carlsen-Larsen-Sims-Vittadello:Co-universal}*{Section~4} and \cite{doron2020cenvelopes}*{Definition~4.6}]\label{def:gauge-comp} Let $\pi=\{\pi_p\}$ be a representation of $\Hilm=(\Hilm_p)_{p\in P}$ in a $\Cst$\nb-algebra $B$ and suppose that $\Cst(\pi)$ admits a coaction $\gamma\colon\Cst(\pi)\to \Cst(\pi)\otimes\Cst(G)$. We will say that~$\pi$ is \emph{gauge-compatible} with~$\gamma$ if  for all $p\in P$ and $\xi\in \Hilm_p$ we have  $$\gamma(\pi_p(\xi))=\pi_p(\xi)\otimes u_p.$$
\end{defn}

Notice that if $(B, G, \gamma)$ is a coaction and $\pi=\{\pi_p\}_{p\in P}$ is a representation of $\Hilm$ in $B$, then $\pi $ is gauge-compatible with $\gamma$ if and only if the induced \Star homomorphism $\tilde{\pi}\colon\Fl\to B$ is $\widetilde{\delta}-\gamma$\nb-equivariant, that is, $$\gamma\circ\tilde{\pi}=(\tilde{\pi}\otimes\id_{\Cst(G)})\circ\widetilde{\delta}.$$ In case $\pi$ is gauge-compatible with $\gamma$ and $B=\Cst(\pi)$, it follows that $(B, G, \gamma)$ is necessarily nondegenerate \cite[Lemma~2.2]{kaliszewski_quigg_2016}. When the coaction $\gamma$ is understood, we will simply say that $\pi$ is gauge-compatible.

\begin{cor}\label{cor:gen-co-universal} Let $(\Toepr(\Hilm)_g)_{g\in G}$ be the Fell bundle associated to the gauge coaction~$\bar{\delta}$ on~$\Toepr(\Hilm)$. Then the $\Cst$\nb-envelope $\Cst_{\mathrm{env}}(\Toepr(\Hilm)^+)$ satisfies the following properties:
\begin{enumerate}\label{cor:co-universal}
\item there is a coaction $\delta_\Lambda\colon\Cst_{\mathrm{env}}(\Toepr(\Hilm)^+)\to \Cst_{\mathrm{env}}(\Toepr(\Hilm)^+)\otimes\Cst(G)$ for which the representation of $\Hilm$ induced by the inclusion 
$\iota\colon\Toepr(\Hilm)^+\hookrightarrow \Cst_{\mathrm{env}}(\Toepr(\Hilm)^+)$ is gauge-compatible;

\item  if  $(B, G, \gamma)$ is a coaction and $\pi=\{\pi_p\}_{p\in P}$ is an injective representation of $\Hilm$ in $B$ that is gauge-compatible with $\gamma$ and induces a surjective \Star homomorphism $\hat{\pi}\colon\Cst((\Toepr(\Hilm)_g)_{g\in G})\to B$, then there exists a $\gamma-\delta_\Lambda$\nb-equivariant surjective \Star homomor\-phism $\rho\colon B\to \Cst_{\mathrm{env}}(\Toepr(\Hilm)^+)$ such that the diagram
 $$\xymatrix{
\Cst((\Toepr(\Hilm)_g)_{g\in G}) \ar@{->}[r]^{\hat{\pi}} \ar@{->}[dr]_{\phi_\lambda\circ\Lambda}&
B \ar@{->}[d]_{\rho}\\
&
 \Cst_{\mathrm{env}}(\Toepr(\Hilm)^+)
    }$$ commutes, where $\Lambda\colon \Cst((\Toepr(\Hilm)_g)_{g\in G})\to \Cst_r((\Toepr(\Hilm)_g)_{g\in G})\cong\Toepr(\Hilm)$ is the left regular representation and $\phi_\lambda\colon \Toepr(\Hilm)\to \Cst_{\mathrm{env}}(\Toepr(\Hilm)^+)$ is the quotient map.
\end{enumerate}
\begin{proof} The normal coaction $(\Cst_{\mathrm{env}}(\Toepr(\Hilm)^+), G, \delta_\Lambda)$ as in the statement was already used in the proof of \thmref{thm:main-result} and follows from the description of $\Cst_{\mathrm{env}}(\Toepr(\Hilm)^+)$ as the reduced cross sectional $\Cst$\nb-algebra of the Fell bundle coming from the coaction~$\delta$ on the covariance algebra of $\Hilm$. In order to show that $\Cst_{\mathrm{env}}(\Toepr(\Hilm)^+)$ also has property (2), let $\hat{\pi}\colon\Cst((\Toepr(\Hilm)_g)_{g\in G})\to B$ be a surjective \Star homomorphism induced by an injective representation $\pi=\{\pi_p\}_{p\in P}$ of $\Hilm$ in $B$ that is gauge-compatible with~$\gamma$. Let $(B_g)_{g\in G}$ be the Fell bundle associated to~$(B, G, \gamma)$ and suppose first that $\gamma$ is normal. Then $B\cong \Cst_r((B_g)_{g\in G})$ with an isomorphism that identifies the corresponding spectral subspaces. It follows that $\hat{\pi}$ factors through $\Toepr(\Hilm)$, giving a surjective \Star homomorphism $\hat{\pi}_*\colon \Toepr(\Hilm)\to B$ such that $\hat{\pi}=\hat{\pi}_*\circ\Lambda$. Also, $\hat{\pi}_*$ is injective on the copy of~$A$ and is $\overline{\delta}-\gamma$-equivariant, and so the (faithful) conditional expectation of $B$ onto $B_e$ is compatible with $E_\lambda$. Hence \corref{cor:gauge-iso} implies that the restriction of $\hat{\pi}_*$ to the tensor algebra $\Toepr(\Hilm)^+$ is completely isometric. It follows that there exists a surjective \Star homomorphism $\rho\colon B\to  \Cst_{\mathrm{env}}(\Toepr(\Hilm)^+)$ such that $\rho\circ\hat{\pi}_*\restriction_{\Toepr(\Hilm)^+}=\iota$. Thus $\phi_\lambda=\rho\circ \hat{\pi}_*$, which yields $\phi_\lambda\circ \Lambda=\rho\circ\hat{\pi}$. This gives (2) when $\gamma$ is a normal coaction. The general case follows from this one because there is a surjective \Star homomorphism $B\to\Cst_r((B_g)_{g\in G})$ that identifies the fibres of $(B_g)_{g\in G}$ by \cite[Theorem~19.5]{Exel:Partial_dynamical} and $\Cst_r((B_g)_{g\in G})$ carries a normal coaction of~$G$ for which the spectral subspace at $g\in G$ is precisely the canonical copy of~$B_g$  (see \cite[Proposition~18.24]{Exel:Partial_dynamical} and  \cite[Proposition~3.4]{doron2020cenvelopes}).
\end{proof}
\end{cor}

If $\Hilm=\CC^P$ is the canonical product system over $P$ with one-dimensional fibres, so that $\Toepr(\Hilm)=\Toepr(P)$, then $\Cst((\Toepr(P)_g)_{g\in G})$ can be canonically identified with the universal Toeplitz algebra $\Toepu(P)$ as introduced in~\cite{laca-sehnem}. Theorem~4.2 of \cite{kakariadis2021boundary} shows that $\partial\Toepr(P)$ is also co-universal for gauge-equivariant nonzero representations of Li's semigroup $\Cst$\nb-algebra $\Cst_s(P)$ \cite[Definition~3.2]{Li:Semigroup_amenability}. If $P$ does not satisfy the independence condition, $\Cst((\Toepr(P)_g)_{g\in G})$ is then a proper quotient of $\Cst_s(P)$ (see \cite[Corollary~3.23]{laca-sehnem}). This implies in particular that $\Cst((\Toepr(\Hilm)_g)_{g\in G})$ is in general not the largest $\Cst$\nb-algebra for which $\Cst_{\mathrm{env}}(\Toepr(\Hilm)^+)$ has the co-universal property with respect to gauge-equivariant representations that are injective on the corresponding copy of~$A$. Under the assumption that $\Hilm$ is faithful, we give next a class of representations for which  $\Cst_{\mathrm{env}}(\Toepr(\Hilm)^+)$ has the co-universal property that is in general much larger than the class of injective gauge-compatible representations that induce a \Star homomorphism of $\Cst((\Toepr(\Hilm)_g)_{g\in G})$. 

\begin{thm}\label{thm:zero-element}  Let $\Hilm=(\Hilm_p)_{p\in P}$ be a faithful product system over $P$ with coefficient $\Cst$\nb-algebra~$A$. Let $(B, G, \gamma)$ be a coaction and let $\pi=\{\pi_p\}_{p\in P}$ be an injective representation of $\Hilm$ in $B$ such that $B=\Cst(\pi)$. Suppose, in addition, that $\pi=\{\pi_p\}_{p\in P}$ is gauge-compatible with $\gamma$  and satisfies the following property: for every neutral word $\alpha=(p_1,p_2,\ldots, p_{2k-1},p_{2k})\in \W(P)$ such that $K(\alpha)=\emptyset$ and every choice of elements $\xi_{p_i}\in \Hilm_{p_i}$ for $i=1,\ldots, 2k$, we have $$\pi_{p_1}(\xi_{p_1})\pi_{p_2}(\xi_{p_2})^*\ldots \pi_{p_{2k-1}}(\xi_{p_{2k-1}})\pi_{p_{2k}}(\xi_{p_{2k}})^*=0.$$ Then there exists a $\gamma-\delta_\Lambda$\nb-equivariant surjective \Star homomorphism $\rho\colon B\to  \Cst_{\mathrm{env}}(\Toepr(\Hilm)^+)$ such that $\rho\circ\pi_p=\iota_p$ for all $p\in P$.
\begin{proof} Consider the \Star homomorphism $\tilde{\pi}\colon \Fl\to B$ obtained by the universal property of $\Fl$. We will show that $\ker\tilde{\pi}\cap \Fl_e$ is contained in $J_e$.  To do so, we begin by observing that a simple application of the $\Cst$\nb-axiom shows that $$\pi_{p_1}(\xi_{p_1})\pi_{p_2}(\xi_{p_2})^*\ldots \pi_{p_{2k-1}}(\xi_{p_{2k-1}})\pi_{p_{2k}}(\xi_{p_{2k}})^*=0$$ if $K(\alpha)=\emptyset$ also when $\alpha$ is not neutral since $K(\alpha)=K(\tilde{\alpha}\alpha)$ and $\tilde{\alpha}\alpha$ is then a neutral word (see, for example,  \cite[Proposition~2.6]{laca-sehnem}). 

We claim that if $\alpha=(p_1,p_2,\ldots, p_{2k-1},p_{2k})$ is a neutral word and $F\subset G$ is a finite set containing the iterated quotient set $Q(\alpha)$ of $\alpha$, then for every $r\in P$ such that $r\not\in K(\alpha)$ and for every $\xi_r\in \Hilm_r I_{r^{-1}(r\vee F)}$, we have that $$\pi_{p_1}(\xi_{p_1})\pi_{p_2}(\xi_{p_2})^*\ldots \pi_{p_{2k-1}}(\xi_{p_{2k-1}})\pi_{p_{2k}}(\xi_{p_{2k}})^*\pi_r(\xi_r)=0.$$ Indeed, suppose that $r\not\in K(\alpha)$. Then there is $g\in Q(\alpha)$ such that $r\not\in gP$. If $rP\cap gP=\emptyset$, then $K(\alpha)\cap rP\subset gP\cap rP=\emptyset$. Setting $\beta\coloneqq (p_1,p_2,\ldots, p_{2k-1},p_{2k}, r,e)$, we see that $K(\beta)=\emptyset$ and this yields $$\pi_{p_1}(\xi_{p_1})\pi_{p_2}(\xi_{p_2})^*\ldots \pi_{p_{2k-1}}(\xi_{p_{2k-1}})\pi_{p_{2k}}(\xi_{p_{2k}})^*\pi_r(\xi_r)=0$$ by the assumption and by nondegeneracy of the right action of $A$ on $\Hilm_r$. Now if $rP\cap gP\neq \emptyset$, we can find $s\in rP\cap gP$ and because $\Hilm$ is faithful and $r\not\in gP$, we obtain from the definition of $I_{r^{-1}(r\vee F)}$ that $$I_{r^{-1}(r\vee F)}\subset\ker\varphi_{r^{-1}s}=\{0\}.$$ Hence $\xi_r=0$, giving $$\pi_{p_1}(\xi_{p_1})\pi_{p_2}(\xi_{p_2})^*\ldots \pi_{p_{2k-1}}(\xi_{p_{2k-1}})\pi_{p_{2k}}(\xi_{p_{2k}})^*\pi_r(\xi_r)=0.$$ This proves the claim.

Next take $b\in\ker\tilde{\pi}\cap  \Fl_e$. In order to show that $b\in J_e$, let $\varepsilon>0$ and let $b'\in  \Fl_e$ with $\|b-b'\|<\frac{\varepsilon}{2}$ and such that $b'=\sum_{\substack{i=1}}^nb_i$, where each $b_i$ is of the form $$b_i=\tilde{t}_{p_1}(\xi_{p_1})\tilde{t}_{p_2}(\xi_{p_2})^*\ldots \tilde{t}_{p_{2k-1}}(\xi_{p_{2k-1}}) \tilde{t}_{p_{2k}}(\xi_{p_{2k}})^*$$ and $\alpha_i\coloneqq (p_1,p_2,\ldots, p_{2k-1}, p_{2k})\in\W(P)$ is a neutral word. Let $F\coloneqq \bigcup_{\substack{i=1}}^nQ(\alpha_i)$ and let $\xi=\sum_{\substack{r\in F'}}\xi_r\in\Hilm_r I_{r^{-1}(r\vee F)}$ with $F'\subset P$ finite. Then since $\tilde{\pi}(b)=0$, we have \begin{equation*}
\begin{aligned}
\|\sum_{\substack{r\in F'}}\pi_r(\xi_r)^*\tilde{\pi}(b')^*\tilde{\pi}(b')\pi_r(\xi_r)\|&=\|\sum_{\substack{r\in F'}}\pi_r(\xi_r)^*\tilde{\pi}(b'-b)^*\tilde{\pi}(b'-b)\pi_r(\xi_r)\|\\& <\frac{\varepsilon^2}{4}\|\sum_{\substack{r\in F'}}\pi_e(\braket{\xi_r}{\xi_r}\|.
\end{aligned}
\end{equation*} Now we apply the usual argument: the left-hand side of the inequality above is precisely $\|t_{F}(b')\xi\|^2$ because $\tilde{\pi}(b_i)\pi_r(\xi_r)=0$ unless $r\in K(\alpha_i)$ and $\pi$ is injective on $A$, while $\|\sum_{\substack{r\in F'}}\pi_e(\braket{\xi_r}{\xi_r}\|=\|\xi\|^2$. Hence $\|q(b')\|= \|b'\|_F<\frac{\varepsilon}{2}$, and thus $\|q(b)\|<\varepsilon$, where $q\colon \Fl\to A\times_{\Hilm}P$ is the quotient map. Since $\varepsilon>0$ is arbitrary, we conclude that $q(b)=0$, that is, $b\in J_e$ as wanted.

Finally, let $(B_g)_{g\in G}$ be the Fell bundle associated to $(B, G,\gamma)$. Because $\ker\tilde{\pi}\cap\Fl_e\subset J_e$, it follows that the restriction of the quotient map $q\colon \Fl\to A\times_{\Hilm}P$ to the spectral subspace $\Fl_g$ at $g\in G$ factors through $\tilde{\pi}$. From this we obtain a morphism of Fell bundles $\phi\colon (B_g)_{g\in G}\to ([A\times_{\Hilm}P]_g)_{g\in G}$ such that $\phi\circ \pi_p=j_p$ for all $p\in P$ because $\tilde{\pi}$ is $\widetilde{\delta}-\gamma$\nb-equivariant. Such a morphism induces a \Star homomorphism $$\hat{\phi}\colon \Cst_r((B_g)_{g\in G})\to\Cst_r\big(([A\times_{\Hilm}P]_g)_{g\in G}\big)\cong \Cst_{\mathrm{env}}(\Toepr(\Hilm)^+)$$ by \cite[Proposition~21.3]{Exel:Partial_dynamical}. Since the collection $\{B_g\}_{g\in G}$ is a topological grading for $B$, it follows from \cite[Theorem~19.5]{Exel:Partial_dynamical} that there exists a \Star homomorphism $\psi\colon B\to \Cst_r((B_g)_{g\in G})$ that identifies the fibres of $(B_g)_{g\in G}$. So setting $\rho\coloneqq \hat{\phi}\circ \psi$ we conclude that $\rho\colon B\to  \Cst_{\mathrm{env}}(\Toepr(\Hilm)^+)$ is a \Star homomorphism satisfying $\rho\circ\pi_p=\iota_p$, and hence is $\gamma-\delta_\Lambda$\nb-equivariant. This completes the proof of the theorem.
\end{proof}
\end{thm}

\begin{rem} Whereas in the proof \corref{cor:co-universal} it was more convenient to apply the defining property of $\Cst_{\mathrm{env}}(\Toepr(\Hilm)^+)$ as the smallest $\Cst$\nb-algebra generated by a completely isometric copy of $\Toepr(\Hilm)^+$, the defining relations of the covariance algebra of $\Hilm$ were a crucial tool in the proof \thmref{thm:zero-element} above. This illustrates the importance of the distinct descriptions of $\Cst_{\mathrm{env}}(\Toepr(\Hilm)^+)$ provided by \thmref{thm:main-result}.
\end{rem}

If $\Hilm=\CC^P$, the relation required in the statement of \thmref{thm:zero-element} corresponds to relation~(T2) of \cite[Definition~3.6]{laca-sehnem}. So \thmref{thm:zero-element} together with \corref{cor:boundary-quotient} imply that $\partial\Toepr(P)$ has the co-universal property for a class of nonzero gauge-compatible isometric representations of~$P$  that is in general strictly larger than the class of representations satisfying the defining relations of Li's semigroup $\Cst$\nb-algebra $\Cst_s(P)$ (see \cite[Proposition~3.22]{laca-sehnem}). Thus the next corollary provides a strengthening of \cite[Theorem~4.2]{kakariadis2021boundary}.

\begin{cor} Let $P$ be a submonoid of a group $G$.  Let $(B, G, \gamma)$ be a coaction and suppose that $B$ is generated as a $\Cst$\nb-algebra by a nonzero isometric representation $w\colon P\to B$ that is gauge-compatible with $\gamma$ and satisfies the following property: for every neutral word $\alpha=(p_1,p_2,\ldots, p_{2k-1},p_{2k})\in \W(P)$ such that $K(\alpha)=\emptyset$, we have $$w_{p_1}w_{p_2}^*\ldots w_{p_{2k-1}}w_{p_{2k}}^*=0.$$ Then there exists a surjective \Star homomorphism $\rho\colon B\to  \partial\Toepr(P)$ such that $\rho\circ w_p=j_p$ for all $p\in P$, where $p\mapsto j_p$ is the canonical isometric representation of $P$ in $\partial\Toepr(P)$.
\end{cor}

\begin{rem} We observe that in general $\Cst_{\mathrm{env}}(\Toepr(\Hilm)^+)$ does not have the co-universal property with respect to gauge-equivariant surjective \Star homomorphisms of Fowler's $\Cst$\nb-algebra $\Fl$ that are injective on $A$. Indeed, if $\Hilm=\CC^P$, then $\Fl$ is the universal $\Cst$\nb-algebra for isometric representations of $P$  and hence has the full and reduced group $\Cst$\nb-algebras $\Cst(G)$ and $\Cst_{\lambda}(G)$ as canonical gauge-equivariant quotients if $P\hookrightarrow G$ generates $G$ as a group. If $P$ is not lef reversible, then the \Star homomorphism from $\Fl$ onto $\partial\Toepr(P)$ induced by $j=\{j_p\}_{p\in P}$ does not factor through the group $\Cst$\nb-algebra because $j_pj_p^*j_qj_q^*=0$ in $\partial\Toepr(P)$ if $pP\cap qP=\emptyset$, while $u_pu_p^*u_qu_q^*=1$ in $\Cst(G)$. See \cite[Example~4.7]{Carlsen-Larsen-Sims-Vittadello:Co-universal}.
\end{rem}

\end{document}